\documentclass[11pt]{article}
\newif\ifanonymous
\anonymousfalse

\usepackage[colorlinks = true,linkcolor = blue,urlcolor  = blue,citecolor = blue,anchorcolor = blue]{hyperref}
\usepackage{mathtools}
\usepackage{subcaption}
\usepackage{amsmath, amssymb, amsthm}
\usepackage{natbib}
\usepackage[margin=1in]{geometry}
\usepackage{parskip}
\usepackage{setspace}
\title{Stable Central Limit Theorems for Discrete-time Lag Martingale Difference Arrays: Applications to Dynamic Causal Inference}
\author{%
  \ifanonymous
    \large Anonymous Authors$^*$
  \else
    \large Walter Dempsey$^{*,1}$ \and
    \large Easton Huch$^{*,2}$
  \fi
}
\date{August 7, 2026}

\renewcommand\vec[1]{\boldsymbol{#1}}
\newcommand\mat[1]{\mathbf{#1}}

\DeclareMathOperator{\Cov}{Cov}
\DeclareMathOperator{\Var}{Var}

\DeclareMathOperator{\Sym}{Sym}
\newcommand{\E}{\mathbb{E}}
\newcommand{\R}{\mathbb{R}}
\newcommand{\N}{\mathbb{N}}
\newcommand{\G}{\mathcal{G}}
\newcommand{\Hcal}{\mathcal{H}}
\newcommand{\A}{\mathcal{A}}
\newcommand{\B}{\mathcal{B}}
\newcommand{\X}{\mathcal{X}}

\newcommand{\Z}{\mathbb{Z}}

\newcommand{\F}{\mathcal{F}}
\newcommand{\I}{\mathcal{I}}
\newcommand{\RR}{\mathcal{R}}
\newcommand{\PP}{\mathcal{P}}
\renewcommand{\L}{\mathcal{L}}

\newcommand{\labeleditem}[1]{%
  \item[$\mathrm{(}$\csname cond#1\endcsname$\mathrm{)}$]%
  \hypertarget{cond:#1}{}%
}
\newcommand{\refitem}[1]{$\mathrm{(}$\hyperlink{cond:#1}{\csname cond#1\endcsname}$\mathrm{)}$}

\newcommand{\pconverge}{\overset{\mathrm{p}}{\to}}
\newcommand{\dconverge}{\overset{\mathrm{d}}{\to}}

\newtheorem{theorem}{Theorem}
\newtheorem{proposition}{Proposition}
\newtheorem{lemma}{Lemma}
\newtheorem{corollary}{Corollary}
\newtheorem{remark}{Remark}
\newtheorem*{externaltheorem}{Theorem}
\theoremstyle{definition}
\newtheorem{assumption}{Assumption}
\newtheorem{definition}{Definition}

\begin{document}

\maketitle

\let\thefootnote\relax
\ifanonymous
\footnotetext{$^{*}$Author names are removed for anonymous peer review.}
\else
\let\thefootnote\relax
\footnotetext{$^{*}$The authors contributed equally to this work.}
\footnotetext{$^{1}$Department of Biostatistics, University of Michigan, Ann Arbor, MI 48109, USA. \texttt{wdem@umich.edu}}
\footnotetext{$^{2}$Johns Hopkins Carey Business School, 555 Pennsylvania Ave NW, Washington, DC 20001, USA. \texttt{ehuch@jhu.edu}}
\fi

\begin{abstract}
\normalsize
Recent work in dynamic causal inference introduced a class of discrete-time stochastic processes that generalize martingale difference sequences and arrays as follows: the random variates in each sequence have expectation zero given certain lagged filtrations but not given the natural filtration.
We formalize this class of stochastic processes and prove stable central limit theorems (CLTs) via martingale--coboundary decomposition, leveraging the classical martingale CLT.
We develop a variety of sufficient conditions, including conditions under which the limiting variance has a simple form that depends on variances and covariances of neighboring variates.
We demonstrate the application of these results to inference for time-averaged treatment effects in switchback designs and present a simulation study supporting their validity.
The CLTs enable various extensions to existing methodology for design-based approaches to dynamic causal inference, including time-lagged effects, random limiting variances, cross-unit dependence, and vector-valued estimands.
\end{abstract}

\vfill


\newpage
\section{Introduction}\label{sec:intro}

Recent work by Bojinov, Rambachan, and Shephard (hereafter, BRS) proposes randomization-based methods for estimating time-averaged effects of treatments on future outcomes \citep{bojinov_time_2019,bojinov_panel_2021}.
As a simple example, consider a setting in which an experimenter sequentially randomizes a single individual to a binary treatment $A_t \in \{0, 1\}$ for $T \in \N$ time points.
Following each treatment, the experimenter observes an outcome $Y_t \in \R$.
The proposed methods enable the experimenter to estimate the time-averaged effect of $A_t$ on $Y_{t+p}$ for $p = 0, 1, 2, \ldots$.
With $p=0$, the result is an estimate of the `contemporaneous' or immediate effect of $A_t$ on $Y_t$.
For $p \geq 1$, the resulting estimate provides a measure of effect persistence over time, with significant effects at larger values of $p$ indicating greater persistence.
Assessing treatment persistence is a well-recognized goal in dynamic causal inference, and other authors have proposed similar definitions of persistence, usually in a superpopulation setting
\citep{RobinsJamesM.1999EotC,BLACKWELLMATTHEW2018HtMC,Boruvka03072018}.

BRS propose asymptotically Gaussian inference procedures justified by the martingale central limit theorem (CLT) \citep{hall_martingale_1980,hausler_stable_2015};
however, the martingale CLT does not apply for $p \geq 1$.
We show in Section \ref{sec:application} that the theory fails for $p \geq 1$ because the supposed-martingale sequence, $(U_t)_{1 \leq t \leq T-p}$, generally does not satisfy the \emph{martingale property}, $\E(U_t | \Hcal_{t-1}) = 0$, where $(\Hcal_t)_{0 \leq t \leq T-p}$ is a filtration with each $U_t$ being $\Hcal_t$-measurable.
Rather, it satisfies a condition that we refer to as the \emph{lag martingale property}: $\E(U_t | \Hcal_{t-p-1}) = 0$.
For $p \geq 1$, we show in \eqref{eq:cov} that the covariance between neighboring time points is nonnegligible; hence, the martingale CLT does not apply.

In spite of this theoretical issue, the proposed framework of BRS poses considerable merit because it produces unbiased estimates of dynamic causal effects under weak assumptions.
In particular, the framework imposes only weak regularity conditions on the data-generating process for the outcomes.
This stands in contrast to many competing methods in dynamic causal inference that require restrictive conditions, such as `no carryover effects' across time periods \citep{ImaiKosuke2019WSWU,deChaisemartinClément2020TFEE,arkhangelsky2022doubly,ATHEY202262}.
In addition, the BRS framework does not require a superpopulation assumption or independent units, enabling its application to a finite set of experimental units.

These benefits have motivated numerous recent proposals applying design-based estimation techniques to problems in dynamic causal inference.
These proposals extend the basic design-based framework of BRS, modifying it to accommodate adaptive treatment assignment \citep{ham2023design, li2026design}, temporal switchback designs \citep{hu2022switchback,bojinov2023design,zeng2026sequentially}, spatiotemporal interference \citep{ni2023design}, and regression-based estimators \citep{lin2025unifying}.
These methods and accompanying experimental designs have been employed at numerous firms, including Lyft \citep{chamandy2016experimentation}, DoorDash \citep{kastelman2018switchback}, and Procter \& Gamble \citep{ni2025reliable}.
Consequently, addressing the above theoretical challenges of the BRS framework could provide significant benefit to methodologists and practitioners, enabling valid asymptotic inference for the case $p \geq 1$ under weak assumptions---most notably, without imposing strong restrictions on carryover effects.

The primary contributions of this work are (i) formulating the class of lag martingales that arise in these problems, (ii) connecting this class to related notions and techniques in probability theory (most notably, stable convergence), (iii) developing corresponding stable CLTs with simple sufficient conditions and an interpretable limiting variance, and (iv) illustrating these results by applying them to the BRS framework in simulation.
Similar to the classical martingale CLT, our results guarantee stable convergence to a scale mixture of Gaussian random variables.
An intuitive characterization of stable convergence is \citep[Proposition VIII.5.33]{alma991060420510207861}
\[
X_n \to X \ \G\text{-stably} \quad \iff \quad (X_n, Y) \dconverge (X, Y) \quad \text{for all } \G \text{-measurable } Y,
\]
where $\dconverge$ denotes (weak) convergence in distribution.
Stability enables us to obtain a limiting Gaussian distribution (as opposed to a scale-mixture of Gaussian distributions) by applying the stable Cram\'er--Slutsky rule \citep[Theorem 3.18 (b)]{hausler_stable_2015}: 
\[
X_n \to X\  \G\text{-stably},\quad Y_n \pconverge Y,\quad Y \text{ is } \G\text{-measurable} \quad \implies \quad X_n Y_n \dconverge X Y
\]
where $\pconverge$ denotes convergence in probability.
In dynamic causal inference, $1/Y^2$ is the variance of the estimator and may be a random variable.
The randomness of this variance arises due to the possibility that early treatments affect future outcome and treatment trajectories.

From a technical standpoint, our theoretical results rely on a martingale--coboundary decomposition.
This technique, attributable to \citet{gordin1969central}, decomposes the partial sums of interest into a sum of martingale differences and an asymptotically negligible coboundary term, enabling the application of classical martingale CLTs.
While initially proposed for stationary processes, the technique was subsequently extended to nonstationary processes \citep{philipp1975almost,volny2006}, such as those studied here.
Whereas classical martingale--coboundary decompositions involve a full telescoping sum, we modify the technique slightly to allow variations between neighboring terms, enabling us to handle varying lag lengths.

An alternative approach is to apply Bernstein blocking, partitioning the rows of a triangular array into minor and major blocks and applying a known CLT to the major-block sums \citep{bernstein_sur_1927,Dehling2002}.
In this manuscript, we employ the martingale coboundary strategy as it results in simpler theoretical statements and sufficient conditions.

In probability theory, CLTs bearing the greatest similarity to those considered here involve mixingale processes, which constrain conditional expectations according to a decaying sequence of constants.
With a fixed lag length, mixingales are easily seen to nest lag martingales as a special case; however, this is not true when the lag length is allowed to diverge.
\citet{mcleish1975invariance, mcleish1977invariance} introduced mixingales and proved CLTs for them, including in the relatively more challenging case of nonstationary processes.
His early results have subsequently been extended and refined by \citet{Wooldridge_White_1988}, \citet{deJong_1997}, \citet{IkedaShinS.2017Anom}, and others.
Recent work by \citet{merlevede2019} provides functional CLTs for nonstationary processes under conditions less restrictive than mixingale conditions.

While the above CLTs can sometimes be applied to lag martingales, our specialized theoretical treatment offers four key benefits.
First, the lag martingale property simplifies the theoretical results and their corresponding conditions.
Second, our approach results in stable CLTs with respect to a $\sigma$-field large enough to include the random limiting variance, which admits a simple and interpretable expression.
Third, our main CLTs (Theorems \ref{thm:basic-clt}--\ref{thm:sufficient-clt}) allow diverging lag lengths, enabling application to more complex estimands or data-collection settings.
Fourth, our results are stated in terms of vector-valued processes, which enables straightforward extension beyond scalar-valued estimands.
In addition to clarifying the theoretical properties of the estimators proposed by BRS, these benefits also enable various extensions of the BRS framework that we plan to address in future work.

The remainder of the paper is organized as follows.
Section \ref{sec:lag-martingale} provides the theoretical setup and definitions of lag martingale difference sequences and arrays.
Section \ref{sec:basic-clt} provides a basic CLT via the classical martingale CLT.
Section \ref{sec:sufficient-conditions} then provides various sufficient conditions for this CLT.
Section \ref{sec:fixed-lag-length} considers the case of a fixed lag length, developing further sufficient conditions for this specialized setting.
Section \ref{sec:application} illustrates the application of our theoretical results to the BRS framework and includes a simulation study supporting their validity.
Section \ref{sec:conclusion} concludes with a brief discussion of the results.

\section{Lag Martingale Difference Sequences and Arrays}\label{sec:lag-martingale}

Let $(k_n)_{n \in \N}$ be a nondecreasing sequence with $k_n \in \N$ and $k_n \geq n$.
Let $(X_{nk})_{1 \leq k \leq k_n}$ be a sequence of real-valued random variables defined on the probability space $(\Omega, \X, P)$ for all $n \in \N$.
Let $(\X_{nk})_{0 \leq k \leq k_n}$ be a filtration in $\X$ so that $\X_{n0} \subset \X_{n1} \subset \ldots \subset \X_{n k_n} \subset \X$.
For all $n \in \N$, we assume that $(X_{nk})_{1 \leq k \leq k_n}$ is adapted to the filtration $(\X_{nk})_{0 \leq k \leq k_n}$, meaning that $X_{nk}$ is $\X_{nk}$-measurable for all $1 \leq k \leq k_n$.
In other words, the triangular array $(X_{nk})_{1 \leq k \leq k_n, n \in \N}$ of random variables is adapted to the triangular array $(\X_{nk})_{0 \leq k \leq k_n, n \in \N}$ of $\sigma$-fields.
We say that $A \in \L^q(P)$ if $\E(|A|^q) < \infty$.
With this notation in place, scalar-valued lag-$p$ martingale difference sequences are defined as follows.
\begin{definition}[Lag-$p$ martingale difference sequence] \label{def:lag-martingale}
Let $p \in \{0, 1, \ldots\}$ and $n \in \N$ such that $p \leq k_n - 1$.
The sequence $(X_{nk})_{p+1 \leq k \leq k_n}$ is a (scalar-valued) \emph{lag-$p$ martingale difference sequence} if $\E(X_{nk} | \X_{n,k-p-1}) = 0$ and $X_{nk} \in \L^1(P)$ for all $p+1 \leq k \leq k_n$.
\end{definition}
Lag-$p$ martingale difference arrays are defined similarly, allowing $n$ to range over $\N$.
Definition \ref{def:lag-martingale} nests martingale difference sequences (and, by extension, arrays) as the special case $p=0$.
We can relate different lag lengths as follows.
\begin{proposition}[Lag relationships] \label{prop:lag-relationships}
Let $p \in \{0, 1, \ldots\}$, $q \in \{p+1, p+2, \ldots\}$, and $n \in \N$ such that $q \leq k_n - 1$.
If $(X_{nk})_{p+1 \leq k \leq k_n}$ is a lag-$p$ martingale, then $(X_{nk})_{q+1 \leq k \leq k_n}$ is a lag-$q$ martingale.
\end{proposition}
\begin{proof}
Apply the tower property of expectations.
\end{proof}
In particular, Proposition \ref{prop:lag-relationships} implies that a martingale difference sequence is a lag-$p$ martingale difference sequence for any $p \in \N$; however, the converse is false in general.

More generally, we can define lag martingales without specifying a fixed lag length.

\begin{definition}[Lag martingale difference sequence] \label{def:lag-p-martingale}
Let \(n \in \N\), and let \(s_n\) satisfy \(1\leq s_n\leq k_n\).
Let $(p_{nk})_{s_n \leq k \leq k_n}$ be a sequence with each $p_{nk} \in \{0, 1, \ldots, k-1\}$.
Then the sequence $(X_{nk})_{s_n \leq k \leq k_n}$ is a (scalar-valued) \emph{lag martingale difference sequence} if \(X_{nk} \in \L^1(P)\) and $\E(X_{nk} | \X_{n,k-p_{nk} - 1}) = 0$ for all $s_n \leq k \leq k_n$.
\end{definition}

As above, we can define \emph{lag martingale difference arrays} by allowing $n$ to range over $\N$.
This definition allows us to include settings in which $p_n \coloneqq \max_{s_n \leq k \leq k_n} p_{nk} \to \infty$ as $n \to \infty$.
We can further generalize these definitions to vector-valued sequences and arrays by replacing the analogous conditions above with $\E(\vec{X}_{nk} | \X_{n,k-p_{nk}-1}) = \vec{0}_d$ (the $d$-dimensional zero vector) and $X_{nkj} \in \L^1(P)$ for all $1 \leq j \leq d$, where $\vec{X}_{nk} \in \R^d$ and $X_{nkj}$ is the $j$-th element of this vector.

\section{Basic CLT}\label{sec:basic-clt}

This section develops a basic stable CLT for general vector-valued lag martingale difference arrays.
Subsequent sections consider sufficient conditions and adaptations to settings with a fixed lag length.

\subsection{Setup}\label{sec:basic-clt-setup}

Our strategy is to decompose $\vec{X}_{nk}$ into three terms: $\vec{X}_{nk} = \vec{D}_{nk} + \vec{U}_{nk} - \vec{W}_{nk}$.
The array $(\vec{D}_{nk})_{n \in \N, s_n \leq k \leq k_n}$ forms a martingale difference array with each term defined as
\[
\vec{D}_{nk} \coloneqq \sum_{i=0}^{q_{nk}} \vec{D}_{n,k+i,k}, \qquad \vec{D}_{n,k+i,k} \coloneqq \vec{E}_{n,k+i,k} - \vec{E}_{n,k+i,k-1}, \qquad \vec{E}_{njl} \coloneqq \E\left(\vec{X}_{nj} | \X_{nl}\right),
\]
where $\vec{X}_{nr}=\vec{0}_d$ for $r>k_n$ and $q_{nk}$ is defined recursively by
\[
q_{nk}
\coloneqq
\min\left\{
q\in\{q_{n,k-1}, q_{n,k-1} + 1,\ldots\}: q \geq p_{n,k+q}
\right\}
\]
with base case $q_{n,s_n-1} = 0$ (for $k > k_n$, we set $p_{nk} = 0$).
By definition, this implies that the final expectation of the above sum is zero: $\vec{E}_{n,k+q_{nk},k-1} = \vec{0}_d$.
The terms $\vec{U}_{nk}$ and $\vec{W}_{nk}$ comprise the coboundary and are defined as
\[
\vec{U}_{nk} \coloneqq \sum_{i=0}^{q_{nk} - 1} \vec{E}_{n,k+i,k-1}, \qquad \vec{W}_{nk} \coloneqq \sum_{i=1}^{q_{nk}} \vec{E}_{n,k+i,k}
\]
so that $\vec{W}_{nk}$ cancels the first term and $\vec{U}_{nk}$, the second in the summation defining $\vec{D}_{nk}$.
Accounting for cancellations and the fact that $\vec{E}_{n,k+q_{nk},k-1} = \vec{0}_d$, the only remaining term is $\vec{E}_{n,k+0,k} = \E(\vec{X}_{nk} | \X_{nk}) = \vec{X}_{nk}$, showing that the decomposition is valid.

We can then decompose the sum $\overline{\vec{X}}_n \coloneqq \sum_{s_n \leq k \leq k_n} \vec{X}_{nk}$ into two terms, $\overline{\vec{D}}_n$ and $\overline{\vec{R}}_n$, where $\overline{\vec{D}}_n \coloneqq \sum_{s_n \leq k \leq k_n} \vec{D}_{nk}$ and $\overline{\vec{R}}_n \coloneqq \sum_{s_n \leq k \leq k_n} \vec{R}_{nk}$ with $\vec{R}_{nk} = \vec{U}_{nk} - \vec{W}_{nk}$.
We will apply a martingale central limit theorem to $\overline{\vec{D}}_n$.
In contrast, we will assume that $\overline{\vec{R}}_n$ is asymptotically negligible as this term consists of differences in neighboring terms.
For a fixed lag length \(p\), the projection horizon is \(q_{nk}=p\) throughout the row.
Then \(\vec{U}_{nk}=\vec{W}_{n,k-1}\) for \(k>s_n\), so the coboundary terms telescope:
\[
\overline{\vec{R}}_n
=
\sum_{k=s_n}^{k_n}(\vec U_{nk}-\vec W_{nk})
=
\vec U_{n,s_n}-\vec W_{n,k_n}
=
\vec U_{n,s_n},
\]
where the last equality follows from the zero extension beyond \(k_n\).
More generally, we have $\vec{U}_{nk} = \vec{W}_{n,k-1} + \sum_{i=q_{n,k-1}}^{q_{nk}-1} \vec{E}_{n,k+i,k-1}$, which provides the representation
\begin{equation}\label{eq:R-bar-remainder}
\overline{\vec{R}}_n
=
\vec{U}_{n,s_n}
+
\sum_{k=s_n+1}^{k_n}
\sum_{i=q_{n,k-1}}^{q_{nk}-1}
\vec{E}_{n,k+i,k-1}.
\end{equation}
The row index in the increment term is \(k+i\): these are precisely the future terms newly admitted when the projection horizon at \(k\) exceeds the horizon at \(k-1\).
The summation contains at most \(p_n\) terms because the projection-horizon increments telescope.

\subsection{Application of Martingale CLT}\label{sec:martingale-clt}

For a symmetric positive semidefinite matrix $\mat{Q}$, we let $\mat{Q}^{1/2}$ denote the (unique) symmetric positive semidefinite square-root matrix such that $\mat{Q}^{1/2} \mat{Q}^{1/2} = \mat{Q}$.
Let $\|\vec{M}\|$ denote the $\ell^2$ norm of the vector $\|\vec{M}\| \in \R^d$.
The array $(\vec{X}_{nk})_{n \in \N, s_n \leq k \leq k_n}$ is square-integrable if $X_{nkj} \in \L^2(P)$ for all $n \in \N$, $s_n \leq k \leq k_n$, and $1 \leq j \leq d$.
Define
\[
\G_n\coloneqq\bigcap_{m\geq n}\X_{m,k_n},
\qquad
\G\coloneqq\sigma\left(\bigcup_{n=1}^{\infty}\G_n\right).
\]
These $\sigma$-fields are well defined because \(k_m\geq k_n\) for \(m\geq n\), and \((\G_n)_n\) is increasing since, within every row, \(\X_{m,k_n}\subseteq\X_{m,k_{n+1}}\).

We use the following stable martingale CLT as the starting point for our specialized CLTs.
\begin{externaltheorem}[\citet{hausler_stable_2024}, Theorem 2.3]\label{thm:external-hl}
Let \((\vec M_{nk})_{1\leq k\leq k_n,n\in\N}\) be a square-integrable \(\R^d\)-valued martingale difference array adapted to \((\mathcal F_{nk})_{0\leq k\leq k_n,n\in\N}\), with \(k_n\) nondecreasing and \(k_n\geq n\).
Define \(\mathcal G_n^M=\bigcap_{m\geq n}\mathcal F_{m,k_n}\) and \(\mathcal G^M=\sigma(\cup_n\mathcal G_n^M)\).
If, for every \(\epsilon>0\),
\[
\sum_{k=1}^{k_n}
\E\left(
\|\vec M_{nk}\|^2
1\{\|\vec M_{nk}\|\geq\epsilon\}
\mid\mathcal F_{n,k-1}
\right)
\pconverge0
\]
and, for some \(\mathcal G^M\)-measurable symmetric positive semidefinite random matrix \(\mat A\),
\[
\sum_{k=1}^{k_n}
\E(\vec M_{nk}\vec M_{nk}^{\top}\mid\mathcal F_{n,k-1})
\pconverge\mat A,
\]
then \(\sum_{k=1}^{k_n}\vec M_{nk}\to\mat A^{1/2}\vec Z_d\) \(\mathcal G^M\)-stably, where \(\vec Z_d\) is standard Gaussian and independent of \(\mathcal G^M\).
\end{externaltheorem}

By the equivalence of norms on $\R^d$, the conditional Lindeberg assumption required by this Theorem can be modified to use another norm.
We now provide a basic CLT that provides conditions under which we can apply the classical martingale central limit theorem to $(\vec{X}_{nk})_{n \in \N, s_n \leq k \leq k_n}$.

\begin{theorem} \label{thm:basic-clt}
Let $(\vec{X}_{nk})_{n \in \N, s_n \leq k \leq k_n}$ be a square-integrable lag martingale difference array adapted to $(\X_{nk})_{0 \leq k \leq k_n, n \in \N}$.
Assume that
\begin{enumerate}
    \labeleditem{CLB} $\sum_{k=s_n}^{k_n} \E\left(\big\|\vec{D}_{nk}\big\|^2 \cdot 1\left\{\big\|\vec{D}_{nk}\big\| \geq \epsilon\right\} \Big| \X_{n,k-1}\right) \pconverge 0$ as $n \to \infty$ for every $\epsilon > 0$,
    \labeleditem{N} $\sum_{k=s_n}^{k_n} \E(\vec{D}_{nk} \vec{D}_{nk}^{\top} | \X_{n,k-1}) \pconverge \mat{\Psi}$ as $n \to \infty$ for some $\G$-measurable symmetric positive semidefinite random matrix $\mat{\Psi}$, and
    \labeleditem{R} $\overline{\vec{R}}_n \pconverge \vec{0}_d$ as $n \to \infty$.
\end{enumerate}
Then $\overline{\vec{X}}_n \to \mat{\Psi}^{1/2} \vec{Z}_d$ $\G$-stably as $n \to \infty$, where $\vec{Z}_d$ is a standard Gaussian random vector that is independent of $\G$.
\end{theorem}

\begin{proof}
Conditions \refitem{CLB} and \refitem{N} are sufficient to apply Theorem \ref{thm:external-hl} to the zero-padded array with \(\vec D_{nk}=\vec0_d\) for \(1\leq k<s_n\).
Thus, $\overline{\vec{D}}_n \to \mat{\Psi}^{1/2} \vec{Z}_d$ $\G$-stably as $n \to \infty$, where $\vec{Z}_d$ is a standard Gaussian random vector that is independent of $\G$.
We have $\overline{\vec{X}}_n = \overline{\vec{D}}_n + \overline{\vec{R}}_n$ with $\overline{\vec{R}}_n \pconverge \vec{0}_d$ as $n \to \infty$, so the conclusion follows from Theorem 3.18 (a) of \citet{hausler_stable_2015}.
\end{proof}

Conditions \refitem{CLB} and \refitem{N} are common conditions in the literature on martingale CLTs.
Condition \refitem{R} is, in essence, a statement that the coboundary terms are asymptotically negligible.
The primary limitation of Theorem \ref{thm:basic-clt} is that these conditions are stated in terms of $\vec{D}_{nk}$ and $\overline{\vec{R}}_n$ rather than the original variates.

\section{Sufficient Conditions}\label{sec:sufficient-conditions}

This section provides sufficient conditions for \refitem{CLB}, \refitem{N}, and \refitem{R} from Theorem \ref{thm:basic-clt}.
The below conditions are more primitive in that they directly involve the original variates ($\vec{X}_{nk}$) and their expectations ($\vec{E}_{n,k+i,l}$).

\subsection{Sufficient Conditions for the Conditional Lindeberg Assumption}\label{sec:clb-sufficient}

Letting $\PP_{nk} \coloneqq \left\{(i, j) \in \Z^2:\ 0 \leq i \leq q_{nk},\ 0 \leq j \leq 1\right\}$, we now provide sufficient conditions for the \refitem{CLB} assumption.

\begin{lemma}\label{lemma:clb-sufficient}
For $\epsilon > 0$, define the following quantities:
\begin{align*}
    C_{n1} \coloneqq &\ \sum_{k=s_n}^{k_n} \E\left[ \left(\sum_{i=0}^{q_{nk}} \big\| \vec{D}_{n,k+i,k} \big\|\right)^2 \cdot 1\left\{\sum_{i=0}^{q_{nk}} \big\| \vec{D}_{n,k+i,k} \big\| \geq \epsilon\right\} \bigg| \X_{n,k-1}\right], \\
    C_{n2} \coloneqq &\ \sum_{k=s_n}^{k_n} \E\left[(q_{nk} + 1)^2 \cdot \max_{0 \leq i \leq q_{nk}} \big\| \vec{D}_{n,k+i,k} \big\|^2 \cdot 1\left\{\max_{0 \leq i \leq q_{nk}} \big\| \vec{D}_{n,k+i,k} \big\| \geq \frac{\epsilon}{q_{nk} + 1}\right\} \Big| \X_{n,k-1}\right], \\
    C_{n3} \coloneqq &\ \sum_{k=s_n}^{k_n} \E\left[4 (q_{nk} \hspace{-2pt} + \hspace{-2pt} 1)^2 \hspace{-2pt} \cdot \hspace{-2pt} \max_{(i,j) \in \PP_{nk}} \hspace{-1pt} \big\| \vec{E}_{n,k+i,k-j}\big\|^2 \hspace{0pt} \cdot \hspace{-1pt} 1\left\{\max_{(i,j) \in \PP_{nk}} \big\|  \vec{E}_{n,k+i,k-j} \big\| \geq \frac{\epsilon}{2\, (q_{nk} \hspace{-2pt} + \hspace{-2pt} 1)}\right\} \hspace{-2pt} \Big| \X_{n,k-1}\right].
\end{align*}
Denote by $(\mathrm{CLB}_j)$ the condition that $C_{nj} \pconverge 0$ as $n \to \infty$ for any $\epsilon > 0$. Then $(\mathrm{CLB}_3) \implies (\mathrm{CLB}_2) \implies (\mathrm{CLB}_1) \implies$ \refitem{CLB}.
\end{lemma}

\begin{proof}
First, we define the following quantity appearing in the statement of \refitem{CLB}:
\[
    C_{n} \coloneqq \sum_{k=s_n}^{k_n} \E\left(\big\|\vec{D}_{nk}\big\|^2 \cdot 1\left\{\big\|\vec{D}_{nk}\big\| \geq \epsilon\right\} \Big| \X_{n,k-1}\right).
\]
We have the following inequalities for the norm in this quantity: 
\begin{equation*}
    \|\vec{D}_{nk}\| = \left\| \sum_{i=0}^{q_{nk}} \vec{D}_{n,k+i,k} \right\| \leq \sum_{i=0}^{q_{nk}} \big\| \vec{D}_{n,k+i,k} \big\| \leq (q_{nk} + 1) \cdot \max_{0 \leq i \leq q_{nk}} \big\| \vec{D}_{n,k+i,k} \big\|.
\end{equation*}
The above is also no greater than
\[
2\, (q_{nk} + 1) \cdot \max_{(i,j) \in \PP_{nk}} \big\|  \vec{E}_{n,k+i,k-j} \big\|.
\]
Now, for $a,b > 0$ with $a \leq b$, we have $a \cdot 1\{a \geq \epsilon\} \leq b \cdot 1\{b \geq \epsilon\}$.
Consequently, we have $C_n \leq C_{n1} \leq C_{n2} \leq C_{n3}$, which implies the conclusion of the lemma.

\end{proof}

\subsection{Sufficient Condition for the Remainder Assumption}\label{sec:remaining-sufficient}

This subsection gives a sufficient condition for \refitem{R} based on the coboundary terms left over after replacing the original array by its martingale projection array.
Let
\[
\RR_n \coloneqq \{(k, i) \in \Z^2:\ s_n \leq k \leq k_n,\ q_{n,k-1} \leq i \leq q_{nk}-1\}
\]
denote the set of indices arising in the remainder term.
Throughout, when an index set is empty, a maximum of nonnegative norms over that set is interpreted as zero.
With this notation in place, we can provide a sufficient condition for \refitem{R} in terms of a maximum over this set.

\begin{lemma}\label{lemma:remainder-sufficient}
    The following condition is sufficient for \refitem{R}:
    \begin{enumerate}
        \labeleditem{Rprime} $p_n \cdot \max_{(k, i) \in \RR_n} \|\vec{E}_{n,k+i,k-1}\| \pconverge 0 \quad$ as $n \to \infty$.
    \end{enumerate}
\end{lemma}

\begin{proof}  
From \eqref{eq:R-bar-remainder}, we have
\[
\|\overline{\vec{R}}_n\|
= \left\|\sum_{(k,i) \in \RR_n} \vec{E}_{n,k+i,k-1}\right\|
\leq p_n \cdot \max_{(k, i) \in \RR_n} \|\vec{E}_{n,k+i,k-1}\| \pconverge 0 \quad\text{as } n \to \infty.
\]
\end{proof}

Condition \refitem{Rprime} ensures that the growth of the individual expectations does not outpace the maximum lag length.
Of course, it would also be sufficient to allow the maximum in \refitem{Rprime} to range over more values of $k$ and $i$.

\subsection{Sufficient Condition for the Nesting Assumption}\label{sec:nesting-sufficient}

For $0 \leq j \leq p_n$ and $s_n\leq k\leq k_n$, define
\[
\I_{nkj}
\coloneqq
\left\{
i\in\Z:
j\leq i\leq k - s_n \text{ and } i \leq q_{n,k-i}
\right\}.
\]
Because $k \geq s_n$ and $q_{nk} \geq 0$, the set \(\I_{nk0}\) contains $0$; it is also bounded above by $k - s_n$, so \(i_{nk}\coloneqq\max(\I_{nk0})\) is well defined.
Because \(q_{n\ell}\) is nondecreasing in \(\ell\), every nonempty
\(\I_{nkj}\) is a consecutive integer block.
If \(i\in\I_{nkj}\) and \(j\leq h\leq i\), then $h \leq i \leq k - s_n$ and $h\leq i\leq q_{n,k-i} \leq q_{n,k-h}$;
thus \(h\in\I_{nkj}\), and therefore
$
\I_{nkj}=\{j,j+1,\ldots,\max(\I_{nkj})\}
$
whenever \(\I_{nkj}\neq\emptyset\).
Moreover, $q_{nk} \leq p_n$ for every $k$: the base case is $q_{n,s_n-1} = 0$, and if $q_{n,k-1} \leq p_n$, then $q = p_n$ is admissible in the definition of $q_{nk}$ because $p_n \geq q_{n,k-1}$ and $p_n \geq p_{n,k+p_n}$.
Consequently, every element of $\I_{nkj}$ is at most $\min(p_n, k - s_n)$; in particular, $i_{nk} \leq \min(p_n, k - s_n)$.
Since only the lower bound in the definition of $\I_{nkj}$ depends on $j$, we have $\I_{nkj} = \{i \in \I_{nk0} : i \geq j\}$; hence $\I_{nkj} = \{j, j+1, \ldots, i_{nk}\}$ if $j \leq i_{nk}$ and $\I_{nkj} = \emptyset$ otherwise.
In particular, $\max(\I_{nkj}) = i_{nk}$ whenever $\I_{nkj} \neq \emptyset$.

The set $\I_{nkj}$ records which covariance terms appear in the variance decomposition
\[
\sum_{k=s_n}^{k_n} \E(\vec{D}_{nk} \vec{D}_{nk}^\top | \X_{n,k-1}).
\]
Renaming the summation index to $m$ and expanding each $\vec{D}_{nm} = \sum_{i'=0}^{q_{nm}} \vec{D}_{n,m+i',m}$ writes this sum as a sum of products of the form $\E(\vec{D}_{n,m+i',m} \vec{D}_{n,m+i'',m}^\top | \X_{n,m-1})$ taken over $s_n \leq m \leq k_n$ and $0 \leq i', i'' \leq q_{nm}$.
When $i' \geq i''$, the substitution $k = m + i'$, $i = i'$, and $j = i' - i''$ turns such a product into
\[
  \mat{\Psi}_{nkij} = \E(\vec{D}_{n,k,k-i} \vec{D}_{n,k-j,k-i}^\top | \X_{n,k-i-1}),
\]
the conditional covariance of the increments that $\vec{X}_{nk}$ and $\vec{X}_{n,k-j}$ receive at time $k-i$; when $i' < i''$, the same substitution with $i'$ and $i''$ exchanged gives $\mat{\Psi}_{nkij}^\top$.
The term $\mat{\Psi}_{nkij}$ therefore appears in the decomposition exactly when three conditions hold:
(a) $k - i \geq s_n$, so that $\vec{D}_{n,k-i}$ is one of the martingale differences in the row sum;
(b) $i \leq q_{n,k-i}$, so that $\vec{D}_{n,k,k-i}$ is one of the terms defining $\vec{D}_{n,k-i}$; and
(c) $i \geq j$, so that $i'' \geq 0$.
The remaining constraints are automatic: $i' \geq 0$ follows from (c), $i'' \leq q_{nm}$ follows from (b) because $i'' \leq i'$, and products with $k > k_n$ vanish by the zero extension beyond $k_n$.
Requiring (a)--(c) simultaneously gives exactly the definition of $\I_{nkj}$; in particular, both upper bounds on $i$ must hold, so $j \leq i \leq \min(q_{n,k-i}, k - s_n)$.
Since $i \geq j$, condition (a) also gives $k - j \geq s_n$, so both $\vec{X}_{nk}$ and $\vec{X}_{n,k-j}$ lie in the array.

Define the matrix $\mat{\Psi}_n$ as
\begin{equation}\label{eq:psi-n}
\begin{aligned}
\mat{\Psi}_n \coloneqq
\sum_{k=s_n}^{k_n}
&\Var(\vec{X}_{nk} | \X_{n,k-i_{nk}-1})\\
&+
2 \sum_{k=s_n}^{k_n}
\sum_{j=1}^{i_{nk}}
\Sym\left\{
\Cov(\vec{X}_{nk}, \vec{X}_{n,k-j} | \X_{n,k-i_{nk}-1})
\right\},
\end{aligned}
\end{equation}
where $\Sym(\mat{A}) = (\mat{A} + \mat{A}^{\top})/2$ for an arbitrary square matrix, $\mat{A}$.
We can then show the following.

\begin{lemma}\label{lemma:nesting-sufficient}
The following condition is sufficient for \refitem{N}:
\begin{enumerate}
    \labeleditem{Nprime} (i) $\mat{\Psi}_n \pconverge \mat{\Psi}\quad$as $n \to \infty$ with $\mat{\Psi}$ a $\G$-measurable symmetric positive semidefinite random matrix and \\
    (ii) $\left\|\sum_{k=s_n}^{k_n} \E\left(\vec{D}_{nk}\, \vec{D}_{nk}^{\top} | \X_{n,k-1}\right) - \mat{\Psi}_n\right\| \pconverge 0 \quad$as $n \to \infty$.
\end{enumerate}
Further, the term $\sum_{k=s_n}^{k_n} \E\left(\vec{D}_{nk}\, \vec{D}_{nk}^{\top} | \X_{n,k-1}\right) - \mat{\Psi}_n$ admits a representation as a summation of mean-zero variates.
\end{lemma}

\begin{proof}
Condition \refitem{N} requires  $\sum_{k=s_n}^{k_n} \E(\vec{D}_{nk} \vec{D}_{nk}^{\top} | \X_{n,k-1}) \pconverge \mat{\Psi}$.
Thus, the statement follows from the decomposition
\[
\sum_{k=s_n}^{k_n} \E(\vec{D}_{nk} \vec{D}_{nk}^{\top} | \X_{n,k-1}) = \mat{\Psi}_n + \left(\sum_{k=s_n}^{k_n} \E\left(\vec{D}_{nk}\, \vec{D}_{nk}^{\top} | \X_{n,k-1}\right) - \mat{\Psi}_n\right),
\]
in which the first term converges in probability to $\mat{\Psi}$ and the second, to $\mat{0}$ as $n \to \infty$ by (i) and (ii), respectively.

To prove the mean-zero statement, we can decompose the original summation as follows:
\begin{align}
\sum_{k=s_n}^{k_n} \E\left(\vec{D}_{nk}\, \vec{D}_{nk}^{\top} | \X_{n,k-1}\right)
=&\ \sum_{k=s_n}^{k_n} \sum_{i=0}^{q_{nk}} \sum_{i'=0}^{q_{nk}} \E\left(\vec{D}_{n,k+i,k}\, \vec{D}_{n,k+i',k}^{\top} | \X_{n,k-1}\right) \nonumber\\
=&\ \sum_{k=s_n}^{k_n} \sum_{i=0}^{i_{nk}} \mat{\Psi}_{nki0}
+ \sum_{k=s_n}^{k_n}
\sum_{j=1}^{i_{nk}}
\sum_{i=j}^{i_{nk}}
\left(\mat{\Psi}_{nkij} + \mat{\Psi}_{nkij}^{\top}\right),\label{eq:nesting-decomp}
\end{align}
where $\mat{\Psi}_{nki0}$ and $\mat{\Psi}_{nkij}$ are variance and covariance terms defined as
\[
\mat{\Psi}_{nkij} \coloneqq \E\left(\vec{D}_{n,k,k-i}\, \vec{D}_{n,k-j,k-i}^{\top} | \X_{n,k-i-1}\right).
\]
We compute an expectation of these terms as follows:
\begin{align*}
&\ \sum_{i=j}^{i_{nk}} \E \left( \mat{\Psi}_{nkij} \Big| \X_{n,k-i_{nk}-1}\right) \\
= &\ \sum_{i=j}^{i_{nk}} \E \left\{ \Cov\left(\vec{E}_{n,k,k-i}, \vec{E}_{n,k-j,k-i} \Big| \X_{n,k-i-1}\right) \Big| \X_{n,k-i_{nk}-1}\right\} \\
= &\ \Cov(\vec{X}_{nk}, \vec{X}_{n,k-j} | \X_{n,k-i_{nk}-1}),
\end{align*}
where last line follows from iterative application of the law of total covariance. Recall that for $\sigma$-fields $\F \subseteq \Hcal$, the law of total covariance states:
\[\Cov(\vec{A}, \vec{B} | \F) = \E\left\{\Cov\left(\vec{A}, \vec{B} \big| \Hcal\right) \Big| \F\right\} + \Cov\left\{(\E(\vec{A} \big| \Hcal), \E(\vec{B} | \Hcal) \Big| \F\right\}.\]
Applying this with $\F = \X_{n,k-i_{nk}-1}$ and $\Hcal = \X_{n,k-i_{nk}}$:
\begin{align*}
&\ \Cov(\vec{X}_{nk}, \vec{X}_{n,k-j} | \X_{n,k-i_{nk}-1}) \\
=&\ \E[\Cov(\vec{X}_{nk}, \vec{X}_{n,k-j} | \X_{n,k-i_{nk}}) | \X_{n,k-i_{nk}-1}] \\
&\ + \Cov(\vec{E}_{n,k,k-i_{nk}}, \vec{E}_{n,k-j,k-i_{nk}} | \X_{n,k-i_{nk}-1}).
\end{align*}
Because the nonempty set \(\I_{nkj}\) is the consecutive block \(\{j,j+1,\ldots,i_{nk}\}\), iterating this expansion on the first term with progressively finer $\sigma$-fields $\X_{n,k-i}$ for $i = i_{nk}-1, \ldots, j$ yields the covariance terms in the sum. The iteration terminates because $\vec{X}_{n,k-j}$ is $\X_{n,k-j}$-measurable, which implies $\Cov(\vec{X}_{nk}, \vec{X}_{n,k-j} | \X_{n,k-j}) = \mat{0}$.
Note that this derivation includes the variance terms as the special case
\[
\sum_{i=0}^{i_{nk}} \E\left(\mat{\Psi}_{nki0} \Big| \X_{n,k-i_{nk}-1}\right) = \Var(\vec{X}_{nk} | \X_{n,k-i_{nk}-1}).
\]
Matching up terms, we then have
\begin{align*}
    &\ \sum_{k=s_n}^{k_n} \E\left(\vec{D}_{nk}\, \vec{D}_{nk}^{\top} | \X_{n,k-1}\right) - \mat{\Psi}_n \\
    =&\ \sum_{k=s_n}^{k_n} \sum_{i=0}^{i_{nk}} \left\{\mat{\Psi}_{nki0} - \E\left(\mat{\Psi}_{nki0} | \X_{n,k-i_{nk}-1}\right) \right\} \\
    &\ \quad + 2 \sum_{k=s_n}^{k_n}
    \sum_{j=1}^{i_{nk}}
    \sum_{i=j}^{i_{nk}}
    \Sym\left\{\mat{\Psi}_{nkij} - \E\left(\mat{\Psi}_{nkij} | \X_{n,k-i_{nk}-1}\right)\right\},
\end{align*}
where each term is easily seen to be mean-zero, proving the statement.
\end{proof}

In our intended application, \refitem{Nprime} is a convenient primitive condition for \refitem{N} because the covariance matrix $\mat{\Psi}_n$ is defined in terms of covariances between subsequent iterates in the original array.
In contrast, \refitem{N} is stated in terms of the projected iterates.
Moreover, part (ii) in \refitem{Nprime} is natural in that it requires a weak law of large numbers for the centered projection-product array displayed above.

\subsection{A Moment Condition for the Bracket Comparison}\label{sec:moment-bracket-comparison}

This section introduces moment bounds that serve as primitive sufficient conditions for part (ii) of \refitem{Nprime} in Lemma \ref{lemma:nesting-sufficient}.
The moment-growth condition we assume also implies \refitem{CLB} and \refitem{R}, and therefore it provides a direct route to the stable CLT.
Let
\[
\mat{H}_n
\coloneqq
\sum_{k=s_n}^{k_n}
\E(\vec{D}_{nk}\vec{D}_{nk}^{\top}\mid\X_{n,k-1})
\]
denote the projected martingale bracket.

\begin{lemma}[Projection-product moment bound]\label{lemma:projection-product}
Let \(\delta\in(0,2]\) and \(\alpha=1+\delta/2\).
For
\[
\mat{Z}_{nuva}
\coloneqq
\E\left(
\vec{D}_{nua}\vec{D}_{nva}^{\top}
\mid \X_{n,a-1}
\right),
\]
there is a constant \(C<\infty\), depending only on \(d\) and \(\delta\), such that
\[
\E\|\mat{Z}_{nuva}\|_F^\alpha
\leq
C\left\{
\E\|\vec{X}_{nu}\|^{2+\delta}
+
\E\|\vec{X}_{nv}\|^{2+\delta}
\right\}.
\]
\end{lemma}

\begin{proof}
By Jensen's inequality,
\[
\|\mat{Z}_{nuva}\|_F^\alpha
\leq
\E\left(
\|\vec{D}_{nua}\|^\alpha
\|\vec{D}_{nva}\|^\alpha
\mid \X_{n,a-1}
\right).
\]
Taking expectations, applying H\"older's inequality, and using \(2\alpha=2+\delta\) gives
\[
\E\|\mat{Z}_{nuva}\|_F^\alpha
\leq
\left(\E\|\vec{D}_{nua}\|^{2+\delta}\right)^{1/2}
\left(\E\|\vec{D}_{nva}\|^{2+\delta}\right)^{1/2}.
\]
Conditional expectation is a contraction in \(L^{2+\delta}\), so
\[
\|\vec{D}_{nua}\|_{2+\delta}
\leq
\|\vec{E}_{nua}\|_{2+\delta}
+
\|\vec{E}_{n,u,a-1}\|_{2+\delta}
\leq
2\|\vec{X}_{nu}\|_{2+\delta},
\]
and similarly for \(\vec{X}_{nv}\).
The result follows from \(xy\leq (x^2+y^2)/2\).
\end{proof}

\begin{lemma}[Small martingale arrays]\label{lemma:small-md}
Let \(\alpha\in(1,2]\).
Let \((\mat{Y}_{nj},\Hcal_{nj})\) be a triangular array of \(d\times d\) matrix-valued martingale differences, with \(d\) fixed.
There is a constant \(C<\infty\), depending only on \(\alpha\) and \(d\), such that
\[
\left\|\sum_j\mat{Y}_{nj}\right\|_{L^\alpha}
\leq
C\left(\sum_j\E\|\mat{Y}_{nj}\|_F^\alpha\right)^{1/\alpha}.
\]
\end{lemma}

\begin{proof}
Write \(Y_{nj}^{ab}\) for the \((a,b)\) entry of \(\mat{Y}_{nj}\).
For each coordinate pair \((a,b)\), \((Y_{nj}^{ab},\Hcal_{nj})\) is a scalar martingale difference array.
The martingale version of the von Bahr--Esseen inequality \citep[Theorem 2.10]{hall_martingale_1980}, based on \citet{vonbahr_esseen_1965}, gives a constant \(K_\alpha<\infty\) such that
\[
\E\left|\sum_jY_{nj}^{ab}\right|^\alpha
\leq
K_\alpha\sum_j\E|Y_{nj}^{ab}|^\alpha.
\]
Because \(\alpha/2\leq1\),
\[
\left\|\sum_j\mat{Y}_{nj}\right\|_F^\alpha
\leq
\sum_{a,b}\left|\sum_jY_{nj}^{ab}\right|^\alpha.
\]
Taking expectations and applying the coordinatewise bound gives
\[
\E\left\|\sum_j\mat{Y}_{nj}\right\|_F^\alpha
\leq
C_{\alpha,d}\sum_j\E\|\mat{Y}_{nj}\|_F^\alpha.
\]
Taking \(1/\alpha\) powers gives the stated bound.
\end{proof}

\begin{proposition}[Regrouped projection-product bound]\label{prop:regrouped-bound}
Let \(\delta\in(0,2]\), put \(\alpha=1+\delta/2\), and suppose that \(\vec{X}_{nk}\in\L^{2+\delta}(P)\) componentwise for every \(s_n\leq k\leq k_n\).
With
\[
A_n\coloneqq\sum_{k=s_n}^{k_n}\E\|\vec{X}_{nk}\|^{2+\delta},
\qquad
\bar p_n\coloneqq p_n+1,
\]
there is a constant \(C<\infty\), depending only on \(d\) and \(\delta\), such that
\[
\|\mat{H}_n-\mat{\Psi}_n\|_{L^\alpha}
\leq
C\bar p_n^3 A_n^{1/\alpha}.
\]
\end{proposition}

\begin{proof}
Write \(\mat{Z}_{nuva}\) as in Lemma \ref{lemma:projection-product}.
The centered-product identity in Lemma \ref{lemma:nesting-sufficient} gives \(\mat{H}_n-\mat{\Psi}_n\) as a sum of terms
\[
\mat{Z}_{n,k,k,k-i}
-\E(\mat{Z}_{n,k,k,k-i}\mid\X_{n,k-i_{nk}-1})
\]
and, for \(j\geq1\),
\[
\mat{Z}_{n,k,k-j,k-i}
-\E(\mat{Z}_{n,k,k-j,k-i}\mid\X_{n,k-i_{nk}-1}),
\]
together with their transposed orientations.
Since each nonempty \(\I_{nkj}\) is \(\{j,j+1,\ldots,i_{nk}\}\), every centered term can be telescoped into increments
\[
\E(\mat{Z}\mid\X_{n,k-m})
-\E(\mat{Z}\mid\X_{n,k-m-1}),
\qquad j\leq i<m\leq i_{nk},
\]
where \(j=0\) covers the diagonal terms.
For each fixed admissible pair or triple of offsets, write \(r=k-m\).
After zero extension, the corresponding terms form a matrix-valued martingale difference array adapted to \((\X_{nr})_r\).
Conditional Jensen's inequality, the elementary bound for a difference of conditional expectations, and Lemma \ref{lemma:projection-product} imply uniformly over all admissible offsets that
\[
\E\|\mat{Y}_{n,k;\cdot}\|_F^\alpha
\leq
C\left\{
\E\|\vec{X}_{nu}\|^{2+\delta}
+
\E\|\vec{X}_{nv}\|^{2+\delta}
\right\}.
\]
For a fixed diagonal group, summing this display over \(k\) is bounded by \(CA_n\).
For a fixed off-diagonal group, \((u,v)=(k,k-j)\) or \((k-j,k)\), so each row index appears at most twice.
Thus every fixed group satisfies
\[
\sum_k\E\|\mat{Y}_{n,k;\cdot}\|_F^\alpha\leq CA_n.
\]
Lemma \ref{lemma:small-md} gives an \(L^\alpha\) bound of order \(CA_n^{1/\alpha}\) for each fixed group.
Every element of \(\I_{nkj}\) is at most \(p_n\), so all offsets lie in \(\{0,1,\ldots,p_n\}\).
There are therefore at most \(C\bar p_n^2\) diagonal pairs and at most \(C\bar p_n^3\) off-diagonal choices of orientation and \((j,i,m)\).
The \(L^\alpha\) triangle inequality yields the displayed bound.
Appendix \ref{sec:projection-bound-details} gives additional details on the regrouping and counting argument.
\end{proof}

\begin{lemma}[Projection horizon and moment bounds]\label{lemma:moment-conditions}
Suppose that, for some \(\delta\in(0,2]\) and \(\alpha=1+\delta/2\),
\[
\bar p_n^{3\alpha}
\sum_{k=s_n}^{k_n}\E\|\vec{X}_{nk}\|^{2+\delta}
\to0,
\qquad
\bar p_n=p_n+1.
\]
Then \(q_{nk}\leq p_n\), condition \refitem{R} holds, condition \refitem{CLB} holds, and part (ii) of \refitem{Nprime} holds.
\end{lemma}

\begin{proof}
The bound \(q_{nk}\leq p_n\) follows from the induction argument in Section \ref{sec:nesting-sufficient}.
Let \(r=2+\delta\) and \(A_n=\sum_{k=s_n}^{k_n}\E\|\vec{X}_{nk}\|^r\).
The set \(\RR_n\) has cardinality at most \(p_n\), and its row indices \(k+i\) are distinct.
Conditional expectation is an \(L^r\) contraction, so for every \(\epsilon>0\),
\[
\Pr\left\{
p_n\max_{(k,i)\in\RR_n}\|\vec{E}_{n,k+i,k-1}\|>\epsilon
\right\}
\leq
\epsilon^{-r}p_n^rA_n
\to0.
\]
Thus \refitem{Rprime}, and hence \refitem{R}, holds.
For Lindeberg, convexity and \(L^r\)-boundedness of martingale projections give
\[
\sum_{k=s_n}^{k_n}\E\|\vec{D}_{nk}\|^r
\leq
C(p_n+1)^rA_n
\to0.
\]
Therefore, for every \(\epsilon>0\),
\[
\E\left[
\sum_{k=s_n}^{k_n}
\E\left(
\|\vec{D}_{nk}\|^2 1\{\|\vec{D}_{nk}\|\geq\epsilon\}
\mid\X_{n,k-1}
\right)
\right]
\leq
\epsilon^{-\delta}\sum_{k=s_n}^{k_n}\E\|\vec{D}_{nk}\|^r
\to0,
\]
which implies \refitem{CLB}.
Finally, Proposition \ref{prop:regrouped-bound} gives
\[
\|\mat{H}_n-\mat{\Psi}_n\|_{L^\alpha}
\leq
C(\bar p_n^{3\alpha}A_n)^{1/\alpha}
\to0,
\]
which proves part (ii) of \refitem{Nprime}.
\end{proof}

\begin{theorem}[Moment-growth sufficient CLT]\label{thm:moment-growth-clt}
Let $(\vec{X}_{nk})_{n \in \N, s_n \leq k \leq k_n}$ be a square-integrable lag martingale difference array adapted to $(\X_{nk})_{0 \leq k \leq k_n, n \in \N}$.
Assume that part (i) of \refitem{Nprime} holds and, for some \(\delta\in(0,2]\) and \(\alpha=1+\delta/2\),
\[
(p_n+1)^{3\alpha}
\sum_{k=s_n}^{k_n}\E\|\vec{X}_{nk}\|^{2+\delta}
\to0.
\]
Then $\overline{\vec{X}}_n \to \mat{\Psi}^{1/2} \vec{Z}_d$ $\G$-stably as $n \to \infty$, where $\vec{Z}_d$ is a standard Gaussian random vector that is independent of $\G$.
\end{theorem}

\begin{proof}
Lemma \ref{lemma:moment-conditions} verifies \refitem{CLB}, \refitem{R}, and part (ii) of \refitem{Nprime}.
Lemma \ref{lemma:nesting-sufficient} verifies \refitem{N}.
Theorem \ref{thm:basic-clt} gives the result.
\end{proof}

\subsection{Combining Sufficient Conditions}\label{sec:sufficient-cor}

The preceding subsections provide two routes to the assumptions of Theorem \ref{thm:basic-clt}.
Section \ref{sec:moment-bracket-comparison} gives a moment-growth condition that controls the projected bracket comparison and implies the \refitem{CLB} and \refitem{R} assumptions.
In contrast, Sections \ref{sec:clb-sufficient}--\ref{sec:nesting-sufficient} provide separate sufficient conditions for \refitem{CLB}, \refitem{R}, and \refitem{N}, respectively.
Combining these sufficient conditions gives the following theorem for general lag lengths.

\begin{theorem}\label{thm:sufficient-clt}
Let $(\vec{X}_{nk})_{n \in \N, s_n \leq k \leq k_n}$ be a square-integrable lag martingale difference array adapted to $(\X_{nk})_{0 \leq k \leq k_n, n \in \N}$.
Assume that \refitem{Nprime}, \refitem{Rprime}, and any of the conditions $(\mathrm{CLB}_1)$, $(\mathrm{CLB}_2)$, or $(\mathrm{CLB}_3)$ from Lemma \ref{lemma:clb-sufficient} hold.
Then $\overline{\vec{X}}_n \to \mat{\Psi}^{1/2} \vec{Z}_d$ $\G$-stably as $n \to \infty$, where $\vec{Z}_d$ is a standard Gaussian random vector that is independent of $\G$.
\end{theorem}

\begin{proof}
Apply Lemmas \ref{lemma:clb-sufficient}--\ref{lemma:remainder-sufficient} to satisfy the conditions of Theorem \ref{thm:basic-clt}.
\end{proof}

Compared to Theorem \ref{thm:basic-clt}, the conditions of Theorem \ref{thm:sufficient-clt} are more primitive and, in practice, may be easier to check.
Theorem \ref{thm:moment-growth-clt} gives a route to the stable CLT through a single moment-growth condition, while Theorem \ref{thm:sufficient-clt} keeps the verification modular for applications in which different assumptions are easier to check separately.
The next section shows how these conditions can further be simplified in settings with a fixed lag length.

\section{Fixed Lag Length}\label{sec:fixed-lag-length}

In this section, we consider a fixed lag length $p \in \{0, 1, \ldots\}$.
The preceding section treats general lag lengths, while this section records the simplifications that occur when the lag length is fixed.
This case is common in applications and includes the setting considered in Section \ref{sec:application}.
We first provide a sufficient condition for \refitem{CLB} based on the following quantities:
\[
    S_n \coloneqq \max_{(i, j) \in \PP} \sum_{k=s_n}^{k_n} \| \vec{E}_{n,k+i,k-j} \|^2, \qquad
    M_n \coloneqq \max_{(i, j) \in \PP} \max_{s_n \leq k \leq k_n} \| \vec{E}_{n,k+i,k-j} \|,
\]
where $\PP \coloneqq \left\{(i, j) \in \Z^2:\ 0 \leq i \leq p,\ 0 \leq j \leq 1\right\}$.
\begin{lemma}\label{lemma:clb-fixed}
With a fixed lag-length $p \geq 0$, the following condition is sufficient for \refitem{CLB}:
\begin{enumerate}
    \labeleditem{CLBF} $(S_n)_n$ is uniformly integrable, and $M_n \pconverge 0$ as $n \to \infty$.
\end{enumerate}
\end{lemma}

\begin{proof}
We have the inequality
\begin{align*}
&\ \sum_{k=s_n}^{k_n}
\E\left(
\|\vec D_{nk}\|^2
1\{\|\vec D_{nk}\|\geq\epsilon\}
\mid\X_{n,k-1}
\right)\\
\leq{}&
2(p+1)
\sum_{(i,j)\in\PP}
\sum_{(i',j')\in\PP}
\sum_{k=s_n}^{k_n}
\E\left(
\|\vec E_{n,k+i,k-j}\|^2
1\left\{
\|\vec E_{n,k+i',k-j'}\|
\geq\frac{\epsilon}{2(p+1)}
\right\}
\mid\X_{n,k-1}
\right).
\end{align*}
which converges in probability to zero for any $\epsilon > 0$ if
\begin{equation}\label{eq:clb-E-fixed-lag}
\sum_{k=s_n}^{k_n} \E \left(\big\| \vec{E}_{n,k+i,k-j}\big\|^2 \cdot 1\left\{\big\|  \vec{E}_{n,k+i',k-j'} \big\| \geq \epsilon \right\} \Big| \X_{n,k-1} \right) \pconverge 0 \quad \text{as } n \to \infty
\end{equation}
for any $\epsilon > 0$ with $(i, j), (i', j') \in \PP$.
To see that \refitem{CLBF} is sufficient for \eqref{eq:clb-E-fixed-lag}, we will show that 
\begin{equation}\label{eq:clb-fixed-lag-l1}
\sum_{k=s_n}^{k_n} \E \left(\big\| \vec{E}_{n,k+i,k-j}\big\|^2 \cdot 1\left\{\big\|  \vec{E}_{n,k+i',k-j'} \big\| \geq \epsilon \right\}\right) \to 0 \quad \text{as } n \to \infty.
\end{equation}
This strategy suffices because convergence in $L^1$ implies convergence in probability.
We have
\[
\sum_{k=s_n}^{k_n} \E \left(\big\| \vec{E}_{n,k+i,k-j}\big\|^2 \cdot 1\left\{\big\|  \vec{E}_{n,k+i',k-j'} \big\| \geq \epsilon \right\}\right)
\leq \E\left( 1\{M_n \geq \epsilon\} \cdot S_n \right).
\]
Now $1\{M_n \geq \epsilon\} \cdot S_n \pconverge 0$.
Uniform integrability of \(S_n\) gives uniform integrability of \(1\{M_n \geq \epsilon\} \cdot S_n\), so Vitali's theorem gives $\E\left( 1\{M_n \geq \epsilon\} \cdot S_n \right) \to 0$, implying \eqref{eq:clb-fixed-lag-l1} and consequently \eqref{eq:clb-E-fixed-lag}.
\end{proof}

The benefit of \refitem{CLBF} compared to \refitem{CLB} is that \refitem{CLBF} is stated in terms of expectations of the original array.

Lemma \ref{lemma:nesting-sufficient} continues to apply in the
fixed-lag setting. Let
$
\ell_{nk}\coloneqq\min(p,k-s_n).
$
Then the matrix \(\mat{\Psi}_n\) has the exact expression
\[
\mat{\Psi}_n
={}
\sum_{k=s_n}^{k_n}
\Var\!\left(
\vec X_{nk}\mid\X_{n,k-\ell_{nk}-1}
\right)
+
2\sum_{j=1}^{p}\sum_{k=s_n+j}^{k_n}
\Sym\!\left\{
\Cov\!\left(
\vec X_{nk},\vec X_{n,k-j}
\mid\X_{n,k-\ell_{nk}-1}
\right)
\right\}.
\]
Indeed, \(q_{nk} = p\) throughout the row, so \(i_{nk} = \min(p, k - s_n) = \ell_{nk}\), and \eqref{eq:psi-n} reduces to the display above.
Since
$
k-\ell_{nk}=\max(s_n,k-p),
$
the conditioning $\sigma$-field may equivalently be written as
\(\X_{n,\max(s_n,k-p)-1}\).

Away from the first \(p\) row positions, \(\ell_{nk}=p\) and the conditioning $\sigma$-field is \(\X_{n,k-p-1}\). This motivates the
boundary-modified target
\[
\overline{\mat{\Psi}}_n
\coloneqq{}
\sum_{k=s_n}^{k_n}
\Var\!\left(
\vec X_{nk}\mid\X_{n,k-p-1}
\right)
+
2\sum_{j=1}^{p}\sum_{k=s_n+j}^{k_n}
\Sym\!\left\{
\Cov\!\left(
\vec X_{nk},\vec X_{n,k-j}
\mid\X_{n,k-p-1}
\right)
\right\}.
\]
The difference between \(\overline{\mat{\Psi}}_n\) and
\(\mat{\Psi}_n\) is confined to the first \(p\) row positions.
Under condition \refitem{CLBF}, the same uniform-integrability
argument used in Lemma \ref{lemma:clb-fixed}, together with the
standard \(L^1\) bounds for differences of conditional variances and
covariances, gives
$
\left\|
\overline{\mat{\Psi}}_n-\mat{\Psi}_n
\right\|\pconverge0.
$
The application in Section \ref{sec:application} and the structured
panel result in Appendix \ref{sec:structured-panel} verify the
corresponding boundary comparison directly.
Alternatively, in applications it is sometimes convenient to set each $(\X_{nj})_{1 \leq j \leq p}$ equal to $\X_{n0}$, in which case $\overline{\mat{\Psi}}_n = \mat{\Psi}_n$.
Related strategies can also be applied in the general lag setting of Section \ref{sec:sufficient-conditions}.

The following lemma shows that we can further simplify condition \refitem{Rprime} with fixed lag lengths.

\begin{lemma}\label{lemma:remainder-fixed}
With a fixed lag-length $p \geq 0$, the following condition is sufficient for \refitem{R}:
\begin{enumerate}
    \labeleditem{RF} $\max_{0 \leq i \leq p-1} \|\vec{E}_{n,s_n+i,s_n-1}\| \pconverge 0 \quad$ as $n \to \infty$.
\end{enumerate}  
\end{lemma}

\begin{proof}
With a fixed lag length, $\overline{\vec{R}}_n = \vec{U}_{n,s_n}$ consists of only the $p$ terms considered by the above maximum.
Because the number of such terms is finite, we have $\|\overline{\vec{R}}_n\| \pconverge 0$ as $n \to \infty$.
\end{proof}

With fixed lag lengths, Lemma \ref{lemma:remainder-fixed} shows that condition \refitem{R} requires only that the above $p$ expectations converge in probability to zero.

Combining the above sufficient conditions provides the following simplified Theorem.

\begin{theorem}\label{thm:fixed-clt}
Let $(\vec{X}_{nk})_{n \in \N, s_n \leq k \leq k_n}$ be a square-integrable lag martingale difference array adapted to $(\X_{nk})_{0 \leq k \leq k_n, n \in \N}$ with fixed lag length $p \geq 0$.
Assume that conditions \refitem{CLBF}, \refitem{Nprime}, and \refitem{RF} hold.
Then $\overline{\vec{X}}_n \to \mat{\Psi}^{1/2} \vec{Z}_d$ $\G$-stably as $n \to \infty$, where $\vec{Z}_d$ is a standard Gaussian random vector that is independent of $\G$.
\end{theorem}

\begin{proof}
Apply Lemmas \ref{lemma:nesting-sufficient}, \ref{lemma:clb-fixed}, \ref{lemma:remainder-fixed} to show that the conditions of Theorem \ref{thm:basic-clt} are satisfied.
\end{proof}

Although originally proposed for settings with a single experimental unit, the BRS framework can also be extended to panel settings with multiple units, each of which may receive multiple treatments over time \citep{bojinov_panel_2021}.
In these settings, it is useful to distinguish the scientific lag---the within-unit time lag defining the causal estimand---and the martingale lag, $p_{nk}$.
These two lags can differ because the two-dimensional panel must be flattened into a one-dimensional sequence.
Appendix \ref{sec:structured-panel} provides a specialized result for this setting that modifies the above fixed-lag arguments, aggregating contributions observed at the same stage.
This approach preserves the scientific lag and allows certain forms of treatment dependence across stages.

\section{Application}\label{sec:application}

We now apply our theoretical results to the BRS framework described in Section \ref{sec:intro}.

\subsection{Setup}\label{sec:application-setup}

We consider a simplified setting with a single experimental unit representing, for example, a single individual.
The individual is assigned a binary treatment $A_t \in \{0, 1\}$ at each time point $t \in \{1, 2, \ldots, T\}$ for $T \in \N$.
From $t = 2$ onward, we posit potential values \citep{rubin_estimating_1974} for each of these treatment variables.
At $t = 2$, we have $A_2(0)$ and $A_2(1)$, representing the treatment that would be assigned under $A_1 = 0$ or $A_1 = 1$, respectively.
The observed treatment is $A_2 = A_2(A_1)$.
For $t \geq 3$, these potential values depend on the full treatment path up to $t-1$, denoted as $\vec{A}_{t-1} \coloneqq (A_1, A_2, \ldots, A_{t-1})$ or $\vec{a}_{t-1}$ for a fixed (nonrandom) realization.
In general, we have potential treatments $(A_t(\vec{a}_{t-1}))_{\vec{a}_{t-1} \in \{0, 1\}^{t-1}}$ and an observed treatment $A_t = A_t(\vec{A}_{t-1})$.

Associated with each time point, we also have real-valued potential outcomes $(Y_t(\vec{a}_{t}))_{\vec{a}_{t} \in \{0, 1\}^{t}}$ and an observed outcome $Y_t = Y_t(\vec{A}_t)$.
We denote the outcomes up to time point $t$ associated with $\vec{a}_t$ as $\vec{Y}_t(\vec{a}_t) \coloneqq (Y_1(a_1), Y_2(\vec{a}_2), \ldots, Y_t(\vec{a}_t))$, and we denote the observed vector of outcomes up to time point $t$ as $\vec{Y}_t \coloneqq \vec{Y}_t(\vec{A}_t)$.
We denote the full set of potential outcomes across all possible treatment paths as $\mat{Y} \coloneqq (Y_t(\vec{a}_{t}))_{\vec{a}_{t} \in \{0, 1\}^{t}, t \in \{1, 2, \ldots, T\}}$.
As in \citet{bojinov_panel_2021} and \citet{bojinov_time_2019}, we assume that the treatments are sequentially randomized.
\begin{assumption}[Sequential randomization]\label{assumption:sequential-randomization}
For all $t$, $a_t$, and $\vec{a}_{t-1}$, we have
\[
\Pr(A_t = a_t | \vec{A}_{t-1} = \vec{a}_{t-1}, \mat{Y})
= \Pr\left\{A_t = a_t | \vec{A}_{t-1} = \vec{a}_{t-1}, \vec{Y}_{t-1}(\vec{a}_{t-1})\right\}.
\]
\end{assumption}
Below, we use colon notation to denote the value of a variable between certain time points; for example, $\vec{A}_{4:6} \coloneqq (A_4, A_5, A_6)$.
With this notation, we denote the assignment probabilities as follows.
\begin{definition}[Assignment probabilities]\label{def:assignment-probs}
Let $p \in \{0, 1, \ldots, T\}$ and $t \in \{1, 2, \ldots, T-p\}$.
For all $\vec{a}_{t:(t+p)}$, define
\[
p_t\left(\vec{a}_{t:(t+p)} ; \vec{A}_{t-1}, \vec{Y}_{t-1}(\vec{a}_{t-1}) \right) \coloneqq \Pr\left(\vec{A}_{t:(t+p)} = \vec{a}_{t:(t+p)} | \vec{A}_{t-1}, \vec{Y}_{t-1}(\vec{a}_{t-1})\right).
\]
\end{definition}
For brevity, we write $P_t\left(\vec{a}_{t:(t+p)}\right) \coloneqq p_t\left(\vec{a}_{t:(t+p)} ; \vec{A}_{t-1}, \vec{Y}_{t-1}(\vec{a}_{t-1}) \right)$.
We further assume that the assignment probabilities for the observed treatment path are known.
\begin{assumption}[Known assignment probabilities]\label{assumption:known-assignment}
For all $p \in \{0, 1, \ldots, T\}$ and $t \in \{1, 2, \ldots, T-p\}$, the quantity $P_t(\vec{A}_{t:(t+p)})$ is known.
\end{assumption}
We also require positivity for the inverse-probability weights used below.
\begin{assumption}[Assignment positivity]\label{assumption:assignment-positivity}
For each finite row, every treatment-path probability appearing in an inverse-probability denominator is almost surely contained in the interval $(0, 1)$.
\end{assumption}
Strict positivity at each finite row does not by itself control the asymptotic size of the inverse-probability weights.
Assignment probabilities may approach zero only at a rate compatible with the moment condition imposed below.
We define a filtration $(\F_t)_{0 \leq t \leq T}$, where $\F_t \coloneqq \sigma(\mat{Y}, \vec{A}_t)$ is the $\sigma$-field generated by $\mat{Y}$ and $\vec{A}_t$.
Because $\mat{Y}$ is used to generate every $\F_t$, the setup and corresponding statistical procedures effectively fix $\mat{Y}$ at its realized value.

\subsection{Causal Estimand}\label{sec:application-estimand}

Except where noted, the remainder of this section focuses on the special case of $p = 1$ and the following lag-$1$ causal estimand.
\begin{definition}[Lag-$1$ causal effect]\label{def:lag1-causal-effect}
For $t \in \{2, 3, \ldots, T\}$, the lag-$1$ causal effect is
\[
\tau_t \coloneqq
    \frac{\left\{Y_t\left(\vec{A}_{t-2}, 1, 1\right) - Y_t\left(\vec{A}_{t-2}, 0, 1\right)\right\}}{2} +
    \frac{\left\{Y_t\left(\vec{A}_{t-2}, 1, 0\right) - Y_t\left(\vec{A}_{t-2}, 0, 0\right)\right\}}{2}.
\]
\end{definition}
The lag-$1$ causal effect, $\tau_t$, quantifies the effect of treatment $A_{t-1}$ on the outcome at the next time point, $Y_t$, averaging over $A_t = 0$ and $A_t = 1$ with equal probability.
Our statistical goal is to estimate an average of these lag-$1$ causal effects over time.
\begin{definition}[Average lag-$1$ causal effect]\label{def:avg-lag1-causal-effect}
The average lag-$1$ causal effect is
\[\tau \coloneqq \frac{1}{T-1} \sum_{t=2}^T \tau_t.\]
\end{definition}
Definitions \ref{def:lag1-causal-effect} and \ref{def:avg-lag1-causal-effect} correspond with Definition 2 in \citet{bojinov_time_2019} with $p=1$ and uniform weights.
The estimand $\tau$ describes the degree to which treatment effects persist across two consecutive time periods, averaged over time.
For this lag-$1$ application, write
\[
\pi_{t-1}(a)\coloneqq\Pr(A_{t-1}=a\mid\F_{t-2}),
\qquad
\rho_t(b\mid a)\coloneqq\Pr(A_t=b\mid\F_{t-2},A_{t-1}=a).
\]
Then \(P_{t-1}(a,b)=\pi_{t-1}(a)\rho_t(b\mid a)\), and \(P_{t-1}(a)=\pi_{t-1}(a)\) when only one treatment value is displayed.
BRS show that $\tau_t$ and $\tau$ can be unbiasedly estimated according to
\begin{equation}\label{eq:tau-hat}
    \hat{\tau}_t \coloneqq \frac{(-1)^{1-A_{t-1}} Y_t}{2\pi_{t-1}(A_{t-1})\rho_t(A_t\mid A_{t-1})},\quad
    \hat{\tau} \coloneqq \frac{1}{T-1} \sum_{t=2}^T \hat{\tau}_t.
\end{equation}
To study the properties of their general estimators for $p \geq 0$, BRS define a quantity $U_t \coloneqq \hat{\tau}_{t+p} - \tau_{t+p}$ and show that $\E(U_t | \F_{t-1}) = 0$;
however, this definition has the disadvantage that $U_t$ is $\F_{t+p}$-measurable but not $\F_{t}$-measurable unless $p=0$ (in Section \ref{sec:intro}, $\Hcal_t = \F_{t+p}$).
Hence, the classical martingale CLT applies only if $p=0$ because $\E(U_t | \F_{t+p-1}) \neq 0$ in general.
For the lag-$1$ estimand, define $W_t \coloneqq \hat{\tau}_t - \tau_t$, which is $\F_t$-measurable and satisfies $\E(W_t | \F_{t-2}) = 0$.
BRS's assumptions further imply that $W_t \in \L^1(P)$; thus, $(W_t)_{2 \leq t \leq T}$ is a lag-$1$ martingale difference sequence adapted to $(\F_t)_{0 \leq t \leq T}$.

To apply Theorem \ref{thm:moment-growth-clt} with \(p=1\), embed $\mat{Y}$ and $\vec{A}_T$ in an asymptotic regime, write \(W_{nt}\coloneqq\hat\tau_{nt}-\tau_{nt}\), and take
\[
X_{nt}=(n-1)^{-1/2}W_{nt},\qquad s_n=2,\qquad k_n=n,\qquad p_{nt}=1,\qquad n \geq 2.
\]
We scale by $(n-1)^{-1/2}$ because each row consists of $(n-1)$ variates, and the limiting distribution is nondegenerate at inverse-square-root scaling.
By construction, we have the lag martingale property: $\E(W_{nt}\mid\F_{n,t-2})=0$ for $2\leq t\leq n$.
We then assume that
\[
(n-1)^{-1-\delta/2}\sum_{t=2}^{n}\E|W_{nt}|^{2+\delta}\to0
\]
for some $\delta\in(0,2]$.
Also assume that
\begin{equation}\label{eq:psi-n-app}
\overline{\psi}_n
\coloneqq
\frac{1}{n-1}\sum_{t=2}^{n}\Var(W_{nt}\mid\F_{n,t-2})
+
\frac{2}{n-1}\sum_{t=3}^{n}\Cov(W_{nt},W_{n,t-1}\mid\F_{n,t-2})
\pconverge\psi
\end{equation}
for some nonnegative $\G$-measurable random variable $\psi$.
These moment conditions are implied, for example, by bounded potential outcomes and assignment probabilities uniformly bounded away from zero and one.
Uniform positivity is stronger than needed: assignment probabilities may approach zero with \(n\), but only when the resulting inverse weights still satisfy the displayed \((2+\delta)\)-moment rate.

The variance estimator $\overline{\psi}_n$ coincides with $\overline{\mat{\Psi}}_n$ from Section \ref{sec:fixed-lag-length}.
The lower-boundary term in \(\mat{\Psi}_n\) differs from $\overline{\psi}_n$ only at \(t=2\):
\[
\mat{\Psi}_n=\overline{\psi}_n+\Delta_n,
\qquad
\Delta_n
\coloneqq
\frac{1}{n-1}\left\{
\Var(W_{n2}\mid\F_{n1})-\Var(W_{n2}\mid\F_{n0})
\right\}.
\]
Let \(r=2+\delta\) and \(b_n=(n-1)^{-r/2}\, \E|W_{n2}|^r\).
The moment condition gives \(b_n\to0\), and Lyapunov's inequality gives
\[
\E|\Delta_n|
\leq
\frac{2}{n-1}\E|W_{n2}|^2
\leq
2b_n^{2/r}
\to0.
\]
Thus \(\mat{\Psi}_n\pconverge\psi\).
The first conclusion of Theorem \ref{thm:moment-growth-clt} gives
\[
(n-1)^{-1/2}\sum_{t=2}^{n}W_{nt}
\to
\sqrt{\psi}\,Z
\qquad \G\text{-stably}.
\]
If \(\psi>0\) almost surely, define
\[
Z_n^{\mathrm{app}}
\coloneqq
\begin{cases}
\displaystyle
\frac{\sqrt{n-1}(\hat\tau_n-\tau_n)}{\sqrt{\overline{\psi}_n}},
&\overline{\psi}_n>0,\\[3mm]
0,&\overline{\psi}_n\leq0.
\end{cases}
\]
Because \(\Pr(\overline{\psi}_n>0)\to1\), the stable Cram\'er--Slutsky theorem yields
\begin{equation}\label{eq:mixing}
Z_n^{\mathrm{app}}\to Z
\qquad \G\text{-mixing},
\end{equation}
where \(Z\) is standard Gaussian.
This result enables asymptotic Gaussian inference for the average lag-$1$ effect even though \(\overline{\psi}_n\) may converge in probability to a nondegenerate random variable.

In this application, $\overline{\psi}_n$ is not computable from the observed data because it depends on unobserved potential outcomes.
In their work, BRS address this challenge by proposing an estimable variance bound.
Because their original bound addresses only the variance terms, we do not apply it here.
Instead, we directly apply the `oracle' variance estimator, $\overline{\psi}_n$, to validate the theoretical results of Sections \ref{sec:basic-clt}--\ref{sec:fixed-lag-length}.
Future methodological work could address this challenge by developing estimable covariance bounds tailored to particular problem settings and causal estimands.

We proceed by deriving the terms needed to compute $\overline{\psi}_n$ in the simulation study.
In the expressions that follow, $p_t(a_t; \vec{A}_{t-2}, a, \mat{Y}) \coloneqq \Pr(A_t=a_t | \vec{A}_{t-2}, A_{t-1}=a, \mat{Y})$ as in Definition \ref{def:assignment-probs}.
The variance $\Var\left( W_t | \F_{t-2}\right)$ is
\begin{align*}
&\frac{P_{t-1}(0) P_{t-1}(1)}{4} \left\{
\frac{\sum_{a=0}^1 Y_t\left(\vec{A}_{t-2}, 1, a\right)}{P_{t-1}(1)} +
\frac{\sum_{a=0}^1 Y_t\left(\vec{A}_{t-2}, 0, a\right)}{P_{t-1}(0)}
\right\}^2 +\nonumber\\
&\sum_{a=0}^1 \left\{
    \frac{Y_t\left(\vec{A}_{t-2}, a, 1\right)}{p_t\left(1; \vec{A}_{t-2}, a, \mat{Y}\right)} -
    \frac{Y_t\left(\vec{A}_{t-2}, a, 0\right)}{p_t\left(0; \vec{A}_{t-2}, a, \mat{Y}\right)}
\right\}^2
\frac{p_t\left(1; \vec{A}_{t-2}, a, \mat{Y}\right) p_t\left(0; \vec{A}_{t-2}, a, \mat{Y}\right)}{4 P_{t-1}(a)}.\nonumber
\end{align*}
The required covariance term, $\Cov\left( W_{t-1}, W_t | \F_{t-2}\right)$, is
\begin{align}
&\frac{(-1)^{1-A_{t-2}}}{4 P_{t-2}\left(A_{t-2}\right)} \sum_{a_{t-1}=0}^1 \frac{(-1)^{1-a_{t-1}} Y_{t-1}(\vec{A}_{t-2}, a_{t-1})}{P_{t-1}\left(a_{t-1}\right)} \sum_{a_t=0}^1 Y_t(\vec{A}_{t-2}, a_{t-1}, a_t) \label{eq:cov}\\
&- \tau_t \frac{(-1)^{1-A_{t-2}}}{2 P_{t-2}(A_{t-2})} \sum_{a=0}^1 Y_{t-1}(\vec{A}_{t-2}, a).\nonumber
\end{align}

\subsection{Simulation Study}\label{sec:simulation}

We now describe a simulation study we performed to evaluate our theoretical results.

\subsubsection{Simulation Setup.}\label{sec:simulation-setup}
We set $T$=100,000 and perform 10,000 independent replications.
We implemented our simulation in Python, coding the main simulation loop in Cython \citep{behnel2010cython} to reduce computation time.
The simulation completes in under one minute on a personal laptop with 48 GB of RAM and 14 CPUs.
We initially set the assignment probabilities to $0.5$.
After observing $A_t = 1$ at 100 distinct time points, we alter the assignment probability as follows to ensure $\psi$ is a nondegenerate random variable.
If the 100th value of $A_t=1$ occurs at an odd time point, we set
\[
    \Pr(A_t = 1 | \vec{A}_{t-1}, \mat{Y})
    = \begin{cases}
        0.1 & \text{if } A_{t-1} = 1,\\
        0.5 & \text{otherwise}.
    \end{cases}
\]
On the other hand, if the 100th value of $A_t=1$ occurs at an even time point, we set $\Pr(A_t = 1 | \vec{A}_{t-1}, \mat{Y}) = 0.5$ uniformly.
We draw the outcomes as follows:
\[
    Y_t(\vec{A}_{t-1}, a) | \vec{A}_{t-1}  \overset{iid}{\sim} \text{Gamma}\left\{\theta(A_{t-1}), 1\right\},
    \quad \theta(a_{t-1}) = \begin{cases}
        2 & \text{if } a_{t-1} = 0,\\
        3 & \text{if } a_{t-1} = 1
    \end{cases}
\]
for $a = 0, 1$.
This sampling scheme results in a causal effect $\tau$ that differs slightly between replications but satisfies $\E(\tau) = 1$.
Within each replication, $\hat{\tau}$ is conditionally unbiased for the sampled value $\tau$.
In the results below, $n$ is fixed at $T$, so we omit $n$ from subscripts for simplicity.

\begin{figure}[t]
     \centering
     \begin{subfigure}[t]{0.32\textwidth}
         \centering
         \includegraphics[width=\textwidth]{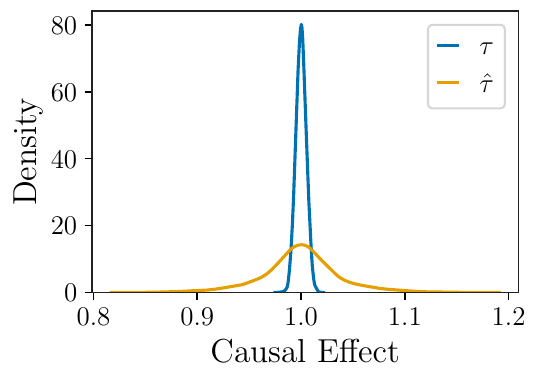}
         \caption{KDEs of $\tau$ and $\Hat{\tau}$}
         \label{fig:tau}
     \end{subfigure}
     \hfill
     \begin{subfigure}[t]{0.32\textwidth}
         \centering
         \includegraphics[width=\textwidth]{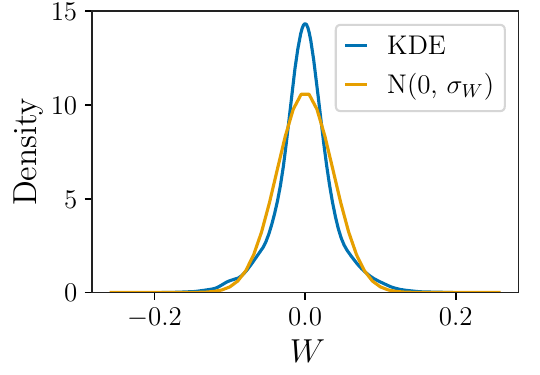}
         \caption{Non-Gaussianity of $W$}
         \label{fig:w}
     \end{subfigure}
     \hfill
     \begin{subfigure}[t]{0.32\textwidth}
         \centering
         \includegraphics[width=\textwidth]{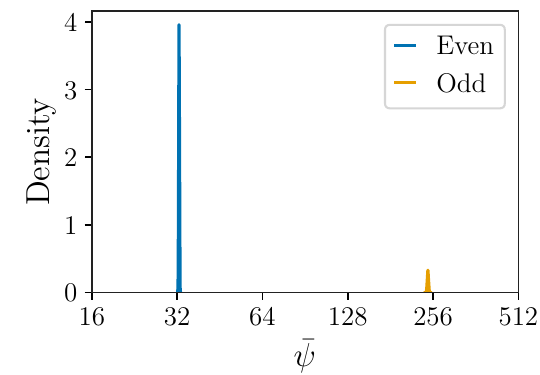}
         \caption{$\overline{\psi}$ by group}
         \label{fig:psi}
     \end{subfigure}
        \caption{Estimates from the simulation study. Panel (a) plots kernel density estimates (KDEs) of the sampled values of $\tau$ and the corresponding unbiased estimator $\Hat{\tau}$. Panel (b) compares a KDE of the sampled values of $W \coloneqq \Hat{\tau} - \tau$ to a Gaussian distribution with variance $\sigma^2_{W} \coloneqq \Var(W)$, offering visual evidence of the non-Gaussianity of $W$. Panel (c) shows the estimates $\overline{\psi}$ segmented by whether the 100th value of $A_t=1$ is even or odd, highlighting how the limiting variance is a random variable.}
        \label{fig:tau-w-psi}
\end{figure}

\subsubsection{Simulation Results.}\label{sec:simulation-results}

Figure \ref{fig:tau-w-psi} displays the estimates from the simulation study.
Panel \ref{fig:tau} plots kernel density estimates (KDEs) of $\tau$ and $\Hat{\tau}$ across replications.
Panel \ref{fig:w} plots a KDE for $W \coloneqq \Hat{\tau} - \tau$ compared to a Gaussian distribution with variance $\sigma^2_W \coloneqq \Var(W)$, providing visual evidence that $W$ is not (marginally) Gaussian.
We also applied the Kolmogorov--Smirnov (KS) test \citep{kolmogorov1933sulla,smirnov1948table} against a Gaussian distribution with variance $\sigma^2_{W}$, resulting in a p-value of $1.5\times 10^{-35}$.

Theorem \ref{thm:moment-growth-clt} implies that $\sqrt{T-1} \, W$ converges stably to a scale mixture of centered Gaussian random variables, with $\psi$ playing the role of random scale parameter.
Panel \ref{fig:psi} plots the simulated values of $\overline{\psi}$ (estimates of $\psi$) segmented by whether the 100th value of $A_t = 1$ occurs at an even or odd time point.
These estimates concentrate around 32.5 and 245.4 for the even and odd groups, respectively, with each group comprising roughly 50\% of the simulation replicates.
Consequently, Theorem \ref{thm:moment-growth-clt} implies that the asymptotic distribution of $\sqrt{T-1} \, W$ is an approximately equal mixture of Gaussian random variables with mean zero and variances corresponding to those of these two groups (estimated at 32.5 and 245.4).
To further evaluate this theoretical implication, we fit a two-component Gaussian mixture model to the sampled values of $\sqrt{T-1} \, W$.
The fitted model estimates one component with a mixture weight of 52\%, mean of $-0.2$, and variance of $241.4$; the other component is estimated at 48\%, -0.1, and $30.8$ for these quantities, respectively.
These estimates closely match the theoretical values we expect based on Theorem \ref{thm:moment-growth-clt}, providing additional empirical support for their validity.

\begin{figure}[t]
     \centering
     \begin{subfigure}[t]{0.49\textwidth}
         \centering
         \includegraphics[width=\textwidth]{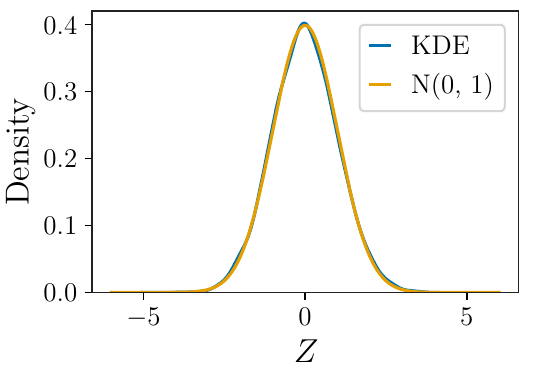}
         \caption{Gaussianity of $Z$}
         \label{fig:z}
     \end{subfigure}
     \hfill
     \begin{subfigure}[t]{0.49\textwidth}
         \centering
         \includegraphics[width=\textwidth]{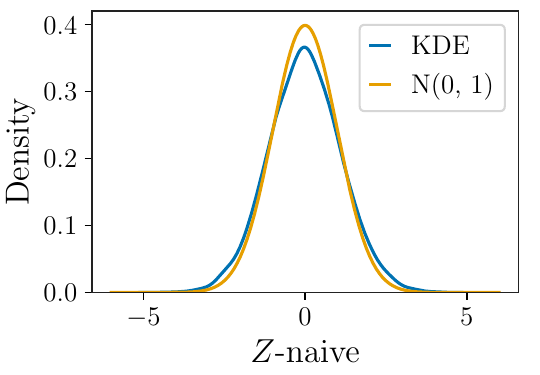}
         \caption{Non-Gaussianity of $Z$-naive}
         \label{fig:z-naive}
     \end{subfigure}
        \caption{Normalized estimates of $W$ from the simulation study. Panel (a) plots a KDE of $Z$ (defined in \eqref{eq:sim-Z}) vs. a standard Gaussian distribution, showing close agreement as expected based on \eqref{eq:mixing}. Panel (b) performs a similar comparison but omitting the covariance terms in the definition of $\overline{\psi}$; we label the resulting quantity $Z$-naive. The comparison of the KDE for $Z$-naive to the standard Gaussian distribution shows significant disagreement, highlighting the need for the covariance terms to obtain valid asymptotic inference.}
        \label{fig:zs}
\end{figure}

To further illustrate the implications of stability in Theorems \ref{thm:basic-clt}--\ref{thm:fixed-clt}, we define
\begin{equation}
\label{eq:sim-Z}
Z \coloneqq \sqrt{\frac{T-1}{\overline{\psi}}} W,
\end{equation}
a realized value from the sequence $\big(Z_n^{\mathrm{app}}\big)_{n \in \N}$.
By \eqref{eq:mixing}, the asymptotic distribution of $Z$ is standard Gaussian.
Figure \ref{fig:zs} plots the empirical distribution of $Z$ and a related quantity, $Z$-naive, that we formed using a variant of $\overline{\psi}$ that omits the covariance terms defined in \eqref{eq:cov}.

Panel \ref{fig:z} shows that the empirical distribution of $Z$ closely matches that of a standard Gaussian as expected.
We applied the KS test and two-sided t-tests for $\E(Z)=0$ and $\E(Z^2)=1$, resulting in p-values of $0.34$, $0.61$, and $0.83$, respectively, offering little to no evidence against the null hypothesis that $Z$ is standard Gaussian.

Panel \ref{fig:z-naive} includes a similar plot for $Z$-naive, showing noticeable disagreement between the empirical distribution and the standard Gaussian distribution.
We applied the hypothesis tests above, resulting in p-values of $0.29$, $0.53$, and $4.2 \times 10^{-37}$, respectively, providing conclusive evidence that the variance is not equal to one.

This comparison highlights the need for the covariance terms in $\overline{\psi}$.
For general $p \geq 1$, their omission results in miscalibrated asymptotic inference.
The resulting confidence intervals may be either conservative or anticonservative depending on the sign of the covariance summation in $\mat{\Psi}_n$.
The naive variance is correct only at $p=0$ because in this case the covariance terms are exactly zero due to the zero-correlation property of martingales.

\section{Conclusion}\label{sec:conclusion}

This paper introduces lag martingales and proves stable CLTs for these processes without requiring restrictive assumptions, such as stationarity.
We provide a simple form for the limiting variance in condition \refitem{Nprime} based on covariances between neighboring variates.

This work was motivated by the dynamic causal inference method of BRS \citep{bojinov_panel_2021,bojinov_time_2019}.
Our theoretical results enable application of their framework with lag length $p \geq 1$.
The stability of our CLTs explicitly allows for a random limiting variance---a possibility that is difficult to rule out in practice.
We illustrate the application of our theoretical results to the BRS framework in Section \ref{sec:application} and provide a simulation study to empirically verify their implications.
Future methodological work in this area could leverage our theoretical results to generalize the BRS framework in various ways, such as allowing cross-unit dependence or extending the methodology to estimate moderation effects.

Although these CLTs are directly applicable to the framework of BRS and its subsequent extensions, we note that the limiting variance typically depends on unobserved potential outcomes.
Thus, in practice, the CLTs do not remove the need to construct estimable bounds of the limiting variance.
They do, however, clarify the form of this limiting variance when estimating lagged treatment effects.

We see four main directions in which our theoretical results could be extended.
First, further sufficient conditions could be developed for specialized settings.
Second, future work could adapt the wide range of conditions for martingale CLTs to provide more flexible limit theorems \citep{hausler_stable_2015,hausler_stable_2024}.
Third, our results could potentially be extended to develop a functional CLT for lag martingales, which could prove useful in developing uniform confidence sequences.
Fourth, future work could extend our setup and results to lag martingales defined in continuous time; this extension could prove useful for data exhibiting continuous fluctuations across time, such as financial time series.

In addition to the technical contributions of this manuscript, our intent is that the application-grounded treatment of the theory provides an accessible introduction to the relevant concepts in stochastic limit theory for researchers who are not experts in this area.
In particular, our simulations in Section \ref{sec:simulation} help elucidate the importance of stability in these results and its connection to the form of the limiting distribution as a Gaussian scale mixture.

\subsection*{Research Data}
The code and results for the simulation found in Section \ref{sec:simulation} are publicly available at 
\ifanonymous
\url{https://github.com/redacted/for/anonymous/peer/review}.
\else
\url{https://github.com/eastonhuch/lag-martingales}.
\fi

\subsection*{Acknowledgements}
\ifanonymous
\else
We thank Gabriel Durham and Inbal Billie Nahum--Shani for helpful feedback and discussions on early versions of this manuscript.
\fi
We thank an anonymous reviewer for suggesting the martingale--coboundary decomposition as an alternative to the Bernstein blocking approach employed in an early version of this manuscript.
The authors acknowledge the use of AI to assist with mathematical derivations, computer code, and writing.
All outputs were manually reviewed for accuracy.

\ifanonymous
\else
\subsection*{Funding}
This work was supported by NIH grants P50DA054039 and R01GM152549 provided by the NIDDK and NIGMS, respectively.
In addition, Dr. Huch gratefully acknowledges research funding from the Johns Hopkins Carey Business School, which supported him during the late stages of this project.
\fi

\bibliographystyle{apalike}
\bibliography{references}

\clearpage
\appendix

\section{Supplement: Details for the Projection-Product Bound}\label{sec:projection-bound-details}

This section provides additional detail for the regrouping step in Proposition \ref{prop:regrouped-bound}.
The proof uses two observations.
First, each centered projection product can be decomposed into martingale differences by telescoping through intermediate \(\sigma\)-fields.
Second, after the offsets are fixed, the resulting terms form a martingale difference array indexed by the remaining row coordinate.

Consider first an off-diagonal term.
For \(j\leq i_{nk}\) and \(i\in\I_{nkj}=\{j,j+1,\ldots,i_{nk}\}\), set
\[
\mat Z_{nkji}
\coloneqq
\E\left(
\vec D_{n,k,k-i}\vec D_{n,k-j,k-i}^{\top}
\mid \X_{n,k-i-1}
\right).
\]
Since \(\mat Z_{nkji}\) is \(\X_{n,k-i-1}\)-measurable and \(i\leq i_{nk}\),
\[
\mat Z_{nkji}
-
\E(\mat Z_{nkji}\mid\X_{n,k-i_{nk}-1})
=
\sum_{m=i+1}^{i_{nk}}
\left\{
\E(\mat Z_{nkji}\mid\X_{n,k-m})
-
\E(\mat Z_{nkji}\mid\X_{n,k-m-1})
\right\}.
\]
The same identity applies to the transposed orientation, and the diagonal terms are the case \(j=0\).

Now fix an admissible choice of offsets \((j,i,m)\) and an orientation.
For example, define
\[
\mat Y_{nk}^{j,i,m}
\coloneqq
\E(\mat Z_{nkji}\mid\X_{n,k-m})
-
\E(\mat Z_{nkji}\mid\X_{n,k-m-1}),
\]
with the convention that \(\mat Y_{nk}^{j,i,m}=\mat 0\) unless \(j\leq i<m\leq i_{nk}\).
Admissibility depends only on \(k\), \(s_n\), and the deterministic sequence \((q_{nk})\), so this convention preserves the martingale difference property below.
As \(k\) varies, the time index \(r=k-m\) varies one-to-one with \(k\).
Thus the array \((\mat Y_{nk}^{j,i,m},\X_{n,k-m})_k\), after zero extension at the row boundaries, is a matrix-valued martingale difference array because
\[
\E(\mat Y_{nk}^{j,i,m}\mid\X_{n,k-m-1})=\mat0.
\]
Conditional Jensen's inequality and the elementary inequality
\[
\| \E(\mat Z\mid\mathcal F)-\E(\mat Z\mid\mathcal F_0)\|_{L^\alpha}
\leq
2\|\mat Z\|_{L^\alpha}
\]
give
\[
\E\|\mat Y_{nk}^{j,i,m}\|_F^\alpha
\leq
C\E\|\mat Z_{nkji}\|_F^\alpha.
\]
Lemma \ref{lemma:projection-product} then yields, uniformly over the fixed offsets,
\[
\E\|\mat Y_{nk}^{j,i,m}\|_F^\alpha
\leq
C\left\{
\E\|\vec X_{nk}\|^{2+\delta}
+
\E\|\vec X_{n,k-j}\|^{2+\delta}
\right\}.
\]
For fixed offsets, summing over \(k\) is therefore bounded by \(CA_n\), since each original row index appears in at most a bounded number of terms.
Lemma \ref{lemma:small-md} gives
\[
\left\|\sum_k\mat Y_{nk}^{j,i,m}\right\|_{L^\alpha}
\leq
CA_n^{1/\alpha}.
\]

It remains only to count the number of fixed-offset groups.
Since \(i_{nk}\leq\min(p_n,k-s_n)\), every admissible offset satisfies \(m\leq i_{nk}\leq p_n\).
The diagonal terms involve pairs \((i,m)\) with \(0\leq i<m\leq p_n\), giving at most \(C(p_n+1)^2\) groups.
The off-diagonal terms involve an orientation, a lag \(1\leq j\leq p_n\), and offsets \(j\leq i<m\leq p_n\), giving at most \(C(p_n+1)^3\) groups.
Applying the \(L^\alpha\) triangle inequality over all groups gives
\[
\|\mat H_n-\mat\Psi_n\|_{L^\alpha}
\leq
C(p_n+1)^3A_n^{1/\alpha},
\]
as claimed in Proposition \ref{prop:regrouped-bound}.

\section{Supplement: A Structured Result for Panel Data}\label{sec:structured-panel}

The main text treats lag martingale difference arrays after they have been indexed by a single row coordinate.
This formulation is appropriate when the lag structure is intrinsic to the one-dimensional array, but panel data introduce a representation issue.
If panel cells are flattened one at a time, a fixed scientific lag can become a row-dependent array lag as the number of individuals grows.
The result below separates this indexing artifact from the substantive lag condition.
By aggregating all contributions observed at the same stage, a fixed scientific lag remains a fixed array lag even when both panel dimensions diverge.
The value of \(p\) is fixed throughout the row limit; the result does not cover \(p=p_n\to\infty\).

Let \(N_n,T_n\to\infty\) be nondecreasing positive integer sequences and put \(S_n=N_n+T_n-1\).
Define
\begin{equation}\label{eq:panel-observation-stage}
\mathcal D_{ns}
=\{(i,t):1\leq i\leq N_n,\ 1\leq t\leq T_n,\ i+t-1=s\},
\qquad 1\leq s\leq S_n.
\end{equation}
Before either side of the rectangle binds, stage \(s\) contains \((1,s),(2,s-1),\ldots,(s,1)\), while in general
\[
|\mathcal D_{ns}|
=
\left[\min(N_n,s)-\max(1,s-T_n+1)+1\right]_+,
\qquad [x]_+=\max(x,0).
\]

For a lag-\(p\) estimand, contributions begin at person-time \(t=p+1\).
Define
\begin{equation}
\A^{(p)}_{ns}
\coloneqq\mathcal D_{ns}\cap\{(i,t):t\geq p+1\},
\qquad
b^{(p)}_{ns}\coloneqq|\A^{(p)}_{ns}|.
\label{eq:panel-contribution-stage}
\end{equation}
\begin{equation}
\begin{aligned}
b^{(p)}_{ns}
&=
\left[\min(N_n,s-p)-\max(1,s-T_n+1)+1\right]_+,\\
M_{np}
&\coloneqq\sum_{s=p+1}^{S_n}b^{(p)}_{ns}
=N_n(T_n-p).
\end{aligned}
\label{eq:panel-block-size}
\end{equation}
Because \(p\) is fixed, \(M_{np}>0\) eventually.
Cells with \(t\leq p\) are observed and may enter later information sets, but do not contribute to the lag-\(p\) sum.

Let \((\Hcal_{ns})_{0\leq s\leq S_n}\) be the stage filtration, with \(\Hcal_{nr}=\Hcal_{n0}\) for \(r\leq0\).
For \((i,t)\in\A^{(p)}_{ns}\), let \(\vec W^{(p)}_{nit}\in\R^d\) be \(\Hcal_{ns}\)-measurable and define
\begin{equation}\label{eq:panel-stage-block}
\vec B^{(p)}_{ns}
=
\sum_{(i,t)\in\A^{(p)}_{ns}}\vec W^{(p)}_{nit},
\qquad p+1\leq s\leq S_n.
\end{equation}
To meet the array convention in Section \ref{sec:lag-martingale} without relabeling rows, put \(L_n\coloneqq\max(S_n,n)\), extend \(\Hcal_{ns}=\Hcal_{n,S_n}\) for \(S_n<s\leq L_n\), and set \(\vec B^{(p)}_{ns}=\vec{0}_d\) for \(S_n<s\leq L_n\).
This zero padding changes neither the panel sum nor its variance target.

The structural condition needed below is the block-level centering condition
\begin{equation}\label{eq:panel-stage-centering}
\E(\vec B^{(p)}_{ns}\mid\Hcal_{n,s-p-1})=\vec{0}_d,
\qquad p+1\leq s\leq S_n.
\end{equation}
A stronger cellwise sufficient condition is
\begin{equation}\label{eq:panel-cellwise-centering}
\E(\vec W^{(p)}_{nit}\mid\Hcal_{n,s-p-1})=\vec{0}_d,
\qquad (i,t)\in\A^{(p)}_{ns}.
\end{equation}
Linearity shows that \eqref{eq:panel-cellwise-centering} implies \eqref{eq:panel-stage-centering}.
Neither condition follows automatically from centering with respect to a person-specific history: enlarging the conditioning \(\sigma\)-field by other individuals' observations can destroy a conditional-mean-zero property.
The global stage-level condition must be justified by the panel design.

\subsection{A Randomization-Based Bridge to Stage Centering}

The next scalar construction gives one route from a global randomization law to \eqref{eq:panel-cellwise-centering}.

\begin{proposition}[Global-history randomization and lag-\(p\) IPW centering]
\label{prop:panel-ipw-centering}
Fix \((i,t)\in\A^{(p)}_{ns}\), so \(s=i+t-1\), and write
\[
\vec A_{n,i,t-p:t}
=(A_{n,i,t-p},A_{n,i,t-p+1},\ldots,A_{nit})\in\{0,1\}^{p+1}.
\]
For every \(\vec a=(a_0,\ldots,a_p)\in\{0,1\}^{p+1}\), suppose there is an integrable path potential outcome \(Y_{nit}(\vec a)\) such that:
\begin{enumerate}
\item \emph{Consistency:}
\[
Y_{nit}^{\mathrm{obs}}=Y_{nit}(\vec a)
\quad\text{on }\{\vec A_{n,i,t-p:t}=\vec a\}.
\]
\item \emph{Measurability:}
For every path, the potential outcome is measurable at the path baseline:
\[
Y_{nit}(\vec a)\text{ is }\Hcal_{n,s-p-1}\text{-measurable}.
\]
\item \emph{Known path randomization and positivity:}
\[
\pi^{(p)}_{nit}(\vec a)
=
\Pr(\vec A_{n,i,t-p:t}=\vec a\mid\Hcal_{n,s-p-1})
\]
is known, \(\Hcal_{n,s-p-1}\)-measurable, and contained in the interval $(0, 1)$ almost surely.
\end{enumerate}
Assume the observed path and \(Y_{nit}^{\mathrm{obs}}\) are \(\Hcal_{ns}\)-measurable.
Define
\begin{align}
\tau^{(p)}_{nit}
&=
\frac1{2^p}
\sum_{(a_1,\ldots,a_p)\in\{0,1\}^p}
\{Y_{nit}(1,a_1,\ldots,a_p)-Y_{nit}(0,a_1,\ldots,a_p)\},
\label{eq:panel-path-estimand}\\
\widehat\tau^{(p)}_{nit}
&=
\frac{(-1)^{1-A_{n,i,t-p}}Y_{nit}^{\mathrm{obs}}}
{2^p\pi^{(p)}_{nit}(\vec A_{n,i,t-p:t})},
\label{eq:panel-path-ipw}\\
W^{(p)}_{nit}
&=\widehat\tau^{(p)}_{nit}-\tau^{(p)}_{nit}.
\end{align}
\vspace{-3pt}Whenever \(\widehat\tau^{(p)}_{nit}\) is integrable, \(W^{(p)}_{nit}\) is \(\Hcal_{ns}\)-measurable and
\vspace{-3pt}
\[
\E(W^{(p)}_{nit}\mid\Hcal_{n,s-p-1})=0.
\]
\end{proposition}

\begin{proof}
Sum over all \(2^{p+1}\) paths and use consistency, measurability, and positivity:
\begin{align*}
\E(\widehat\tau^{(p)}_{nit}\mid\Hcal_{n,s-p-1})
&=
\sum_{\vec a\in\{0,1\}^{p+1}}
\pi^{(p)}_{nit}(\vec a)
\frac{(-1)^{1-a_0}Y_{nit}(\vec a)}
{2^p\pi^{(p)}_{nit}(\vec a)}\\
&=
\frac1{2^p}
\sum_{(a_1,\ldots,a_p)\in\{0,1\}^p}
\Big\{Y_{nit}(1,a_1,\ldots,a_p)
-Y_{nit}(0,a_1,\ldots,a_p)\Big\}\\
&=\tau^{(p)}_{nit}.
\end{align*}
\vspace{-3pt}The centering follows, and adaptation and consistency give measurability.
\end{proof}

The probability in Proposition \ref{prop:panel-ipw-centering} is the
conditional marginal path probability,\vspace{-3pt}
\[
\Pr\!\left(
\vec A_{n,i,t-p:t}=\vec a
\mid\Hcal_{n,s-p-1}
\right),
\]
\vspace{-3pt}induced by the full global assignment mechanism. It is
\(\Hcal_{n,s-p-1}\)-measurable and may depend only on information
available before stage \(s-p\). The global assignment mechanism may be adaptive and dependent across individuals after
that baseline. The proposition does not require this probability to
equal the product of person-specific one-step probabilities, nor
does its centering conclusion require independence of assignments
across individuals, provided the resulting pre-path conditional
marginal probabilities are known and lie in the interval $(0, 1)$. The potential-outcome measurability condition is the randomization-based requirement that the relevant path schedule is fixed before the path begins. 

\subsection{Stage-Blocked Stable Limit}
\label{sec:stage-blocked-stable-limit}

Set
\[
\vec X^{(p)}_{ns}=M_{np}^{-1/2}\vec B^{(p)}_{ns},
\qquad p+1\leq s\leq S_n.
\]
Under \eqref{eq:panel-stage-centering}, this is a lag-\(p\) martingale difference array indexed by stage.
Define
\[
\G_n^{\mathrm{st}}
=\bigcap_{m\geq n}\Hcal_{m,L_n},
\qquad
\G^{\mathrm{st}}
=\sigma\!\left(\bigcup_{n=1}^{\infty}\G_n^{\mathrm{st}}\right).
\]
These \(\sigma\)-fields are well defined because \(L_m\geq L_n\) for \(m\geq n\).
The natural stage-block variance target is
\begin{equation}\label{eq:panel-block-variance}
\begin{aligned}
\mat V_{np}^{\mathrm{st}}
=\frac1{M_{np}}\Bigg[
&\sum_{s=p+1}^{S_n}
\Var(\vec B^{(p)}_{ns}\mid\Hcal_{n,s-p-1})\\
&+
2\sum_{j=1}^{p}\sum_{s=p+1+j}^{S_n}
\Sym\!\left\{
\Cov(\vec B^{(p)}_{ns},\vec B^{(p)}_{n,s-j}
\mid\Hcal_{n,s-p-1})
\right\}
\Bigg].
\end{aligned}
\end{equation}

This target imposes no pairwise independence requirement among cell
contributions within the same stage or between stages separated by
\(j\in\{1,\ldots,p\}\). It includes every corresponding conditional
covariance term, including covariances between contributions from
different individuals, as the following expansions show. This does not
permit unrestricted aggregate dependence: the stage blocks must still
satisfy the centering condition
\eqref{eq:panel-stage-centering}, the block-moment condition
\eqref{eq:panel-block-moment}, and the variance-convergence condition
\eqref{eq:panel-block-variance-convergence}. Blocks separated by more
than \(p\) stages may remain statistically dependent, but the centering
condition makes them uncorrelated.
Writing \(\vec W^{(p)}_{na}=\vec W^{(p)}_{nit}\) for \(a=(i,t)\),
\begin{align*}
\Var(\vec B^{(p)}_{ns}\mid\Hcal_{n,s-p-1})
&=
\sum_{a,a'\in\A^{(p)}_{ns}}
\Cov(\vec W^{(p)}_{na},\vec W^{(p)}_{na'}
\mid\Hcal_{n,s-p-1}),\\
\Cov(\vec B^{(p)}_{ns},\vec B^{(p)}_{n,s-j}
\mid\Hcal_{n,s-p-1})
&=
\sum_{\substack{a\in\A^{(p)}_{ns}\\
a'\in\A^{(p)}_{n,s-j}}}
\Cov(\vec W^{(p)}_{na},\vec W^{(p)}_{na'}
\mid\Hcal_{n,s-p-1}).
\end{align*}
No covariance terms at stage separation greater than \(p\) need to be
added to \eqref{eq:panel-block-variance}. For \(h\geq p+1\), whenever
both blocks are defined, \(\vec B^{(p)}_{n,s-h}\) is
\(\Hcal_{n,s-p-1}\)-measurable. Consequently,
\[
\Cov\!\left(
\vec B^{(p)}_{ns},\vec B^{(p)}_{n,s-h}
\mid\Hcal_{n,s-p-1}
\right)
=\mat{0}_{d\times d}.
\]
Moreover, condition \eqref{eq:panel-stage-centering} implies
\begin{equation}\label{eq:panel-distant-covariance}
\Cov(\vec B^{(p)}_{ns},\vec B^{(p)}_{n,s-h})
=\mat{0}_{d\times d},
\end{equation}
whenever the blocks are square integrable.
Thus distant blocks need not be independent, but they are uncorrelated.

\begin{theorem}[Stage-blocked stable CLT for every fixed finite lag]
\label{thm:panel-stage-block-clt}
Fix \(p\geq0\).
Suppose \eqref{eq:panel-stage-centering} holds and, for some \(\delta\in(0,2]\),
\begin{equation}\label{eq:panel-block-moment}
M_{np}^{-(1+\delta/2)}
\sum_{s=p+1}^{S_n}\E\|\vec B^{(p)}_{ns}\|^{2+\delta}
\longrightarrow0.
\end{equation}
If, for a \(\G^{\mathrm{st}}\)-measurable symmetric positive semidefinite random matrix \(\mat\Psi_p\),
\begin{equation}\label{eq:panel-block-variance-convergence}
\mat V_{np}^{\mathrm{st}}\pconverge\mat\Psi_p,
\end{equation}
then
\[
M_{np}^{-1/2}\sum_{s=p+1}^{S_n}\vec B^{(p)}_{ns}
=
M_{np}^{-1/2}\sum_{i=1}^{N_n}\sum_{t=p+1}^{T_n}\vec W^{(p)}_{nit}
\longrightarrow
\mat\Psi_p^{1/2}\vec Z_d
\qquad \G^{\mathrm{st}}\text{-stably},
\]
where \(\vec Z_d\) is standard Gaussian and independent of \(\G^{\mathrm{st}}\).

In the scalar case, if \(\Psi_p>0\) almost surely, define
\[
T_{n,p}^{\mathrm{st}}
=
\begin{cases}
\displaystyle
\frac{M_{np}^{-1/2}
\sum_{i=1}^{N_n}\sum_{t=p+1}^{T_n}W^{(p)}_{nit}}
{\sqrt{V_{np}^{\mathrm{st}}}},
&V_{np}^{\mathrm{st}}>0,\\[3mm]
0,&V_{np}^{\mathrm{st}}\leq0.
\end{cases}
\]
Then \(T_{n,p}^{\mathrm{st}}\to Z\) \(\G^{\mathrm{st}}\)-mixing, where \(Z\) is standard Gaussian.
\end{theorem}

\begin{proof}
Since \(p\) is fixed and \(N_n,T_n\to\infty\), discard finitely many rows so that \(S_n\geq2p+1\).
Apply Theorem \ref{thm:moment-growth-clt} to the stage-blocked array \(\vec X^{(p)}_{ns}=M_{np}^{-1/2}\vec B^{(p)}_{ns}\), with row beginning \(s_n=p+1\), padded row ending \(L_n\), filtration \(\Hcal_{ns}\), and \(p_{ns}=p\).
Let \(\rho=2+\delta\) and \(\alpha=1+\delta/2\).
Condition \eqref{eq:panel-block-moment} says
\[
A_{np}
\coloneqq
\sum_{s=p+1}^{L_n}\E\|\vec X^{(p)}_{ns}\|^\rho
\longrightarrow0.
\]
Since \((p+1)^{3\alpha}\) is fixed, the moment-growth condition in Theorem \ref{thm:moment-growth-clt} holds.

It remains to verify part (i) of \refitem{Nprime} for the stage-blocked array.
Let \(\mat{\Psi}_{np}^{\mathrm{gen}}\) denote the matrix in \eqref{eq:psi-n} for this array.
Using the fixed-lag form of \eqref{eq:psi-n} and omitting the zero terms from padded rows,
\begin{align*}
\mat{\Psi}_{np}^{\mathrm{gen}}
=&
\sum_{s=p+1}^{S_n}
\Var(\vec X^{(p)}_{ns}
\mid\Hcal_{n,s-\min(p,s-p-1)-1})\\
&+
2\sum_{j=1}^{p}\sum_{s=p+1+j}^{S_n}
\Sym\!\left\{
\Cov(\vec X^{(p)}_{ns},\vec X^{(p)}_{n,s-j}
\mid\Hcal_{n,s-\min(p,s-p-1)-1})
\right\}.
\end{align*}
For \(s=p+1+\ell\), \(0\leq\ell<p\), the generic target conditions on \(\Hcal_{np}\), whereas \eqref{eq:panel-block-variance} conditions on \(\Hcal_{n\ell}\).
For \(\ell\geq p\), the two targets agree.
Thus
\begin{align}
\mat{\Psi}_{np}^{\mathrm{gen}}-\mat V_{np}^{\mathrm{st}}
=&
\sum_{\ell=0}^{p-1}
\{\Var(\vec X^{(p)}_{n,p+1+\ell}\mid\Hcal_{np})
-\Var(\vec X^{(p)}_{n,p+1+\ell}\mid\Hcal_{n\ell})\}
\nonumber\\
&+
2\sum_{\ell=1}^{p-1}\sum_{j=1}^{\ell}
\Sym\Big\{
\Cov(\vec X^{(p)}_{n,p+1+\ell},
\vec X^{(p)}_{n,p+1+\ell-j}\mid\Hcal_{np})
\nonumber\\
&\hspace{38mm}
-\Cov(\vec X^{(p)}_{n,p+1+\ell},
\vec X^{(p)}_{n,p+1+\ell-j}\mid\Hcal_{n\ell})
\Big\}.
\label{eq:panel-exact-boundary-gap}
\end{align}
There is no upper-boundary discrepancy because zero extension excludes out-of-row pairs from both targets.

For square-integrable vectors and \(\sigma\)-fields \(\B_1,\B_2\),
\begin{align*}
\E\|\Var(\vec X\mid\B_1)-\Var(\vec X\mid\B_2)\|_F
&\leq4\E\|\vec X\|^2,\\
\E\|\Cov(\vec U,\vec V\mid\B_1)
-\Cov(\vec U,\vec V\mid\B_2)\|_F
&\leq
4\{\E\|\vec U\|^2\E\|\vec V\|^2\}^{1/2}.
\end{align*}
The second bound follows by expanding conditional covariance and applying conditional Jensen and Cauchy--Schwarz.
Applying these inequalities to \eqref{eq:panel-exact-boundary-gap}, including the outer factor two, gives
\[
\E\|\mat{\Psi}_{np}^{\mathrm{gen}}-\mat V_{np}^{\mathrm{st}}\|_F
\leq
C_p\sum_{\ell=0}^{p-1}\E\|\vec X^{(p)}_{n,p+1+\ell}\|^2
\leq
C_p A_{np}^{2/\rho}
\longrightarrow0.
\]
The last inequality uses Lyapunov's inequality and fixed \(p\).

Equation \eqref{eq:panel-block-variance-convergence} and the \(L^1\) comparison imply \(\mat{\Psi}_{np}^{\mathrm{gen}}\pconverge\mat\Psi_p\), giving part (i) of \refitem{Nprime}.
Theorem \ref{thm:moment-growth-clt} gives the stable CLT.
For scalars, \(V_{np}^{\mathrm{st}}\pconverge\Psi_p>0\) implies \(\Pr(V_{np}^{\mathrm{st}}>0)\to1\).
Changing \(T_{n,p}^{\mathrm{st}}\) on the event \(\{V_{np}^{\mathrm{st}}\leq0\}\) therefore does not affect its limit.
On the event \(\{V_{np}^{\mathrm{st}}>0\}\), the stable Cram\'er--Slutsky theorem gives \(T_{n,p}^{\mathrm{st}}\to Z\) stably.
Since \(Z\) is independent of \(\G^{\mathrm{st}}\), this convergence is \(\G^{\mathrm{st}}\)-mixing.
\end{proof}

\subsection{A Primitive Block-Moment Condition}

\begin{corollary}[A sufficient panel block-moment condition]
\label{cor:panel-block-moment}
Fix \(p\geq0\).
If, for the same \(\delta\) as in Theorem \ref{thm:panel-stage-block-clt}, some \(C<\infty\) satisfies
\begin{equation}\label{eq:panel-rosenthal-bound}
\E\|\vec B^{(p)}_{ns}\|^{2+\delta}
\leq C\{b^{(p)}_{ns}\}^{1+\delta/2}
\qquad\text{for all }n,s,
\end{equation}
then \eqref{eq:panel-block-moment} holds.
\end{corollary}

\begin{proof}
Let \(b_{np}^{\max}=\max_s b^{(p)}_{ns}\).
Since \(\sum_s b^{(p)}_{ns}=M_{np}\),
\[
M_{np}^{-(1+\delta/2)}
\sum_s\E\|\vec B^{(p)}_{ns}\|^{2+\delta}
\leq
C\left(\frac{b_{np}^{\max}}{M_{np}}\right)^{\delta/2}.
\]
For the full rectangle,
\[
b_{np}^{\max}\leq\min(N_n,T_n-p),
\qquad M_{np}=N_n(T_n-p),
\]
so
\[
\frac{b_{np}^{\max}}{M_{np}}
\leq\frac1{\max(N_n,T_n-p)}\longrightarrow0.
\]
\end{proof}

Condition \eqref{eq:panel-rosenthal-bound} follows, for example, from a conditional Rosenthal inequality if, given \(\Hcal_{n,s-p-1}\), the within-stage cell contributions are conditionally independent, conditionally centered, and their conditional \((2+\delta)\)-moments are almost surely bounded by a common constant.
Conditional independence is not necessary: any weak-dependence or martingale ordering that yields the same bound suffices.
Proposition \ref{prop:panel-ipw-centering} does not itself supply this bound.

\begin{remark}[The purely triangular observation schedule]
\label{rem:panel-triangular-schedule}
Suppose observation stops at \(K_n\), where \(N_n\wedge T_n\geq K_n\), \(K_n\) is nondecreasing, and \(K_n-p\to\infty\).
Then
\[
b^{(p)}_{ns}=s-p,\qquad p+1\leq s\leq K_n,
\]
\[
M_{np}=\frac{(K_n-p)(K_n-p+1)}2,
\qquad b_{np}^{\max}=K_n-p,
\]
\[
\frac{b_{np}^{\max}}{M_{np}}
=\frac2{K_n-p+1}\longrightarrow0.
\]
Replace \(S_n\) by \(K_n\), put \(L_n^{K}=\max(K_n,n)\), and zero-pad the blocks and filtration through \(L_n^{K}\).
Use \(\G_n^{\mathrm{st}}=\bigcap_{m\geq n}\Hcal_{m,L_n^{K}}\).
The proof of Theorem \ref{thm:panel-stage-block-clt} is unchanged.
\end{remark}

\begin{remark}[No panel-balance requirement and no growing \(p_n\)]
\label{rem:panel-no-balance}
No relative rate assumption on \(N_n\) and \(T_n\) is needed.
Arbitrarily imbalanced rectangles are allowed because, for fixed \(p\), the geometric moment bound depends only on \(\max(N_n,T_n-p)\to\infty\).
The nondegenerate variance limit remains design-specific.
If \(p=p_n\to\infty\), the generic lag penalty, the number of boundary terms, and the eligible-cell geometry all vary with \(n\).
That regime requires a separate result and is not covered by Theorem \ref{thm:panel-stage-block-clt}.
\end{remark}

\end{document}